\documentclass[a4paper,12pt]{article}

\usepackage[usenames,dvipsnames,svgnames,table]{xcolor}
\usepackage[pagebackref,linktocpage=true,colorlinks=true,linkcolor=Blue,citecolor=BrickRed,urlcolor=RoyalBlue]{hyperref}
\usepackage[alphabetic]{amsrefs}
\usepackage{skak} 
\usepackage{xcolor}

\usepackage{amsmath}
\usepackage{amssymb}
\usepackage{amsthm}
\usepackage{mathtools}
\usepackage{graphicx}
\usepackage{float}
\usepackage{enumitem}

\newtheorem{theorem}{Theorem}[section]

\theoremstyle{remark}
\newtheorem{remark}[theorem]{Remark}

\theoremstyle{definition}

\numberwithin{equation}{section}

\def\R{\mathbb{R}}

\def\d{\partial}
\def\curl{\textnormal{\textrm{curl}}}

\def\dif{{\mathrm d}}
\def\ep{\varepsilon}

\def\x{\vec{x}}
\def\y{\vec{y}}

\def\X{\vec{X}}
\def\vv{\vec{\mathrm{v}}}

\def\hz{\hat{z}}

\renewcommand{\div}{\textnormal{\textrm{div}}}

\begin{document}

\title{
	On explicit geometric solutions of some inviscid flows with free boundary.%
}
\author{
	Yucong {\sc Huang} and
	Aram {\sc Karakhanyan}
}

\date{}
\maketitle

\begin{abstract}
The aim of this paper is to discuss  some well known explicit  examples of the three dimensional inviscid flows with free boundary constructed by John \cite{John} and Ovsyannikov \cite{Ovsyannikov}, and provide a detailed  analysis of its long time behaviour. 
	\end{abstract}

{\hypersetup{linkcolor=blue}
	\tableofcontents
}

\section{Introduction}

The aim of this paper is to discuss  some well known explicit  examples of the three dimensional inviscid flows with free boundary constructed by John \cite{John} and Ovsyannikov \cite{Ovsyannikov}, and provide a detailed  analysis of its long time behaviour. 
It is rare to have such explicit examples of flows even in the two dimensional case \cite{Higgins}, \cite{Gilbarg}, which among other things, makes these solutions special as they describe the instability of the free boundary as $t\to \infty$.
It is important to note  that in all these examples the free boundary, regarded as a surface,  is a conic section. 

In section \ref{sec:2} we study an irrotational inviscid  flow from \cite{Ovsyannikov} where the free boundary is a round cylinder. This example provides an explicit asymptotics in time exhibiting blow up in finite time and shrinking to a line in infinite time.

In section \ref{sec:john} and \ref{sec:4} we provide a similar analysis for the solutions of John \cite{John}. The time asymptotics results are new and did not appear in \cite{John}.

Finally, in section \ref{sec:5} we discuss the case of rotational flow \cite{Ovsyannikov}, which in finite time transform from unit sphere to an ellipsoid and then back to sphere, after which the flow continues as a needle that collapses to a plane at  infinity. 

For each $t\ge 0$, let $\Omega(t)\subset \R^3$ be the moving domain containing the fluid. Denote the Cartesian coordinates as $\x\equiv(x,y,z)$, and let $\vv(t,\x)\vcentcolon= (\mathrm{v}^x,\mathrm{v}^y,\mathrm{v}^z)(t,\x)$ be the velocity vector field. Suppose  the flow is \textbf{incompressible}, i.e. $\div \vv = \d_x \mathrm{v}^x + \d_y \mathrm{v}^y + \d_z \mathrm{v}^z =0$. Moreover, we also assume the flow is \textbf{irrotational}, i.e. $\curl \vv =0$. Then there exists $\phi(t,\x)\vcentcolon [0,\infty)\times\R^3 \to \R$ such that $\vv=\nabla \phi \equiv (\d_x \phi, \d_y \phi, \d_z \phi)$. This then implies that
\begin{equation}\label{Lapla}
	-\Delta \phi \equiv - \big( \d_x^2 \phi + \d_y^2 \phi + \d_z^2 \phi \big) = 0 \qquad \text{for } \ \x\in \Omega(t) \ \text{ and } \ t\ge 0.
\end{equation}
Moreover, the flow at free boundary $\d\Omega(t)$ is modelled by the two postulations:
\begin{enumerate}[label=\textnormal{(\roman*)}]
	\item The free boundary is following the particle path (Kinematic equation);
	\item The Bernoulli's principle is satisfied at the free boundary (Dynamic equation). 
\end{enumerate}  
\paragraph{Formulation of kinematic equation.} Let $t\mapsto \X(t)$ be a point on $\d\Omega(t)$ that moves along with the free boundary. Then postulation \textnormal{(i)} indicates that
\begin{equation}\label{char}
	\dfrac{\dif \X}{\dif t}(t) = \vv(t,\X(t)) = \nabla\phi (t,\X(t)) \quad \text{for } \ t\ge 0 \ \text{ and } \ \X(0)\in \d\Omega(0).
\end{equation}
Suppose that the moving domain is described by the implicit inequality of the form:
\begin{equation*}
	\Omega(t)= \big\{\x\in\R^3  \,\big\vert\, F(t,\x)<0\big\}, \quad \text{at time } \ t\ge 0,
\end{equation*}
for some function $F\vcentcolon\R^{1+3}\to \R$. Then the free boundary is determined by 
\begin{equation*}
\d\Omega(t)= \big\{\x\in\R^3  \,\big\vert\, F(t,\x)=0\big\}, \quad \text{at time } \ t\ge 0.
\end{equation*}
Since $\X(t)$ is a point on the surface, it follows that $F(t,\X(t))=0$ for all $t\ge 0$. Thus by chain rule, we obtain:
\begin{equation}\label{dtF}
	\Big\{ \d_t F  + \nabla \phi \cdot \nabla F \Big\}(t,\X(t)) = 0, \quad \text{for } \ t\ge 0.
\end{equation}
In particular, for each given point $(\x,t)\in \Gamma \vcentcolon= \{ (\tau,\y)\in[0,\infty)\times\R^3  \,\vert\, \y\in \d\Omega(\tau) \}$, let $\{\X(s;\x,t)\}_{0\le s \le t} \in \Gamma$ be the backward characteristic satisfying
\begin{equation}\label{backChar}
	\X(t;t,\x)=\x \quad \text{ and } \quad \dfrac{\d \X}{\d s}(s;t,\x) = \nabla \phi \big( s, \X(s;t,\x) \big), \quad \text{for } \ s\in[0,t]. 
\end{equation}
Then the following differential equation also holds:
\begin{equation}\label{dtFback}
\Big\{ \d_t F + \nabla \phi\cdot \nabla F \Big\}(s,\X(s;t,\x)) = 0, \quad \text{for } \ s\in[0,t].
\end{equation}
\paragraph{Formulation of dynamic equation.} The Bernoulli's principle at the surface is:
\begin{equation*}
	\Big\{\d_t \phi + \dfrac{1}{2}|\nabla \phi|^2 + P \Big\}\Big\vert_{\x\in \d\Omega(t)} = 0,
\end{equation*}
where $P$ is the pressure. In the present paper, we consider the model for which there is no surface tension. Then $P=-c(t)$ for some function of time $c(t)\vcentcolon [0,\infty)\to \R$.  

In summary the incompressible irrotational flow with free boundary is governed by the following system of equations
\begin{subequations}\label{PDEs}
\begin{align}
& - \Delta \phi = 0 && \text{in the domain } \ \x\in \Omega(t).\label{harmonic}\\
& \d_t F + \nabla \phi \cdot \nabla F = 0 && \text{at the surface } \ \x\in\d\Omega(t), \label{lagrange}\\
& \d_t \phi + \dfrac{1}{2} |\nabla \phi|^2 = c(t) \qquad && \text{at the surface } \ \x\in\d\Omega(t).\label{bernoulli}
\end{align}
\end{subequations}
Since the domain $\Omega(t)=\{\x\in\R^3\,\vert\, F(t,\x)< 0\}$ is entirely determined by $F$, (\ref{PDEs}) is a system of $3$ equations with $3$ unknowns, which are $F(t,\x)$, $\phi(t,\x)$, and $c(t)$.

\section*{Acknowledgements}
The authors were partially supported 
by EPSRC grant
EP/S03157X/1 {\em Mean curvature measure of free boundary}.

\section{Flow of Cylinder}\label{sec:2}
In this section, we present the explicit solution to (\ref{PDEs}) first obtained by Ovsjannikov in \cite{Ovsyannikov}. For a function $\alpha(t)$, and constant $\gamma\in\R\backslash\{0\}$, we consider the following ansatz for potential function:
\begin{equation}\label{Ovsjan}
	\phi(t,x,y,z) \vcentcolon= \dfrac{\alpha^{\prime}(t)}{4\gamma\alpha(t)} (x^2+y^2-2z^2).
\end{equation}
It can be easily verified that $\Delta \phi = \d_x^2\phi+\d_y^2\phi +\d_z^2\phi = 0$. Here we are looking for a free boundary taking the form of cylinder, which means that the surface is described by the equation of the form:
\begin{equation}\label{ovsCyl}
	0=F(t,x,y,z)=f\big(t,\sqrt{x^2+y^2}\big) \qquad \text{for some function } f(t,r)\vcentcolon [0,\infty)^2 \to \R.
\end{equation}
\paragraph{Inferences from the surface evolution equation (\ref{lagrange}).} Suppose that the equation $F(t,\x)=0$ describes the free boundary surface. For a given time $t>0$, fix $\x\in \d\Omega(t)$. Let $\{\X(s;t,\x)=(X^x,X^y,X^z)(s;t,\x)\}_{0\le s \le t}$ be the backward characteristic curve emanating from $(t,\x)$, i.e. $s\mapsto \X$ solves (\ref{backChar}). Using the ansatz (\ref{Ovsjan}), one has
\begin{equation*}
	\dfrac{\d X^x}{\d s} = \dfrac{\alpha^{\prime}(s)}{2\gamma\alpha(s)}X^x, \qquad \dfrac{\d X^y}{\d s} = \dfrac{\alpha^{\prime}(s)}{2\gamma\alpha(s)}X^{y}, \qquad \dfrac{\d X^z}{\d s} = - \dfrac{ \alpha^{\prime}(s)}{\gamma\alpha(s)}X^z.
\end{equation*}
Solving the above ODEs we have for $\alpha_0\vcentcolon=\alpha(0)$,
\begin{equation*}
	x \Big(\dfrac{\alpha_0}{\alpha(t)}\Big)^{\frac{1}{2\gamma}} = X^x(0;t,\x), \quad y \Big(\dfrac{\alpha_0}{\alpha(t)}\Big)^{\frac{1}{2\gamma}} = X^y(0;t,\x), \quad z \Big(\dfrac{\alpha(t)}{\alpha_0}\Big)^{\frac{1}{\gamma}}= X^z(0;t,\x).
\end{equation*}
Suppose that at $s=0$, the initial configuration for the free boundary is determined by a given function $\tilde{f}(\x) \in \R^3\to \R$, i.e. $F(0,\x)=\tilde{f}(\x)=0$, Then it follows that 
\begin{equation}\label{gS}
	0=F(t,\x)= \tilde{f}\big(\X(0;t,\x)\big) = \tilde{f}\Big(x \big(\dfrac{\alpha_0}{\alpha(t)}\big)^{\frac{1}{2\gamma}}, y\big(\dfrac{\alpha_0}{\alpha(t)}\big)^{\frac{1}{2\gamma}}, z \big(\dfrac{\alpha(t)}{\alpha_0}\big)^{\frac{1}{\gamma}}\Big).
\end{equation}
Thus the general equations for the free surface takes the form (\ref{gS}), with $\x\to \tilde{f}(\x) $ being a function independent of the variables $(t,x,y,z)$.
\paragraph{Inferences from the Bernoulli's principle (\ref{bernoulli}).}  Substituting the ansatz (\ref{Ovsjan}) into the Bernoulli's principle (\ref{bernoulli}), one has
\begin{equation}\label{gB}
	\Big\{ \dfrac{2\gamma \alpha^{\prime\prime}\alpha + (1-2\gamma)|\alpha^{\prime}|^2 }{8\gamma^2\alpha^2} (x^2+y^2) + \dfrac{(\gamma+1)|\alpha^{\prime}|^2-\gamma\alpha^{\prime\prime}\alpha}{2\gamma^2\alpha^2} z^2 \Big\} \Big\vert_{ F(t,x,y,z)=0} = c(t).
\end{equation}
Since we are seeking free boundaries taking the form of a cylinder described in (\ref{ovsCyl}), we impose that the coefficient for $z^2$ in (\ref{gB}) must be zero, which implies
\begin{equation*}
(\gamma+1)|\alpha^{\prime}|^2 - \gamma \alpha^{\prime\prime}\alpha = 0 \qquad \text{ for } \ t\in\R.
\end{equation*}
Hence $\frac{\alpha^{\prime\prime}}{\alpha^{\prime}}=\frac{\gamma+1}{\gamma}\cdot\frac{\alpha^{\prime}}{\alpha}$, which then gives $\frac{\alpha^{\prime}}{\alpha_0^{\prime}}=(\frac{\alpha}{\alpha_0})^{(\gamma+1)/\gamma}$ where we denote $\alpha_0^{\prime}\vcentcolon=\alpha^{\prime}(0)$. Solving this ODE, it has the solution:
\begin{equation}\label{ovsA}
\alpha(t) = \alpha_0 (1-k_0 t)^{-\gamma} \qquad \text{where } \ k_0 \vcentcolon=\frac{\alpha_0^{\prime}}{\gamma\alpha_0}.
\end{equation} 
Substituting this into (\ref{gB}), one has
\begin{equation}\label{gBCyl}
	\dfrac{3k_0^2}{8(1-k_0 t)^2} (x^2+y^2) \Big\vert_{F(t,x,y,z)=0} = c(t).
\end{equation}
In addition, putting (\ref{ovsA}) into (\ref{gS}), it implies that the free surface must take the form:
\begin{equation}\label{gSCyl}
0 = \tilde{f}\Big( x\sqrt{1-k_0t}, y\sqrt{1-k_0 t}, \dfrac{z}{1-k_0 t} \Big) \quad \text{for some function } \ \tilde{f}\vcentcolon\R^3\to \R.
\end{equation}
Since (\ref{gBCyl}) and (\ref{gSCyl}) must be consistent, it follows that there exists some constant $K\in\R$ such that $c(t) = K (1-k_0 t)^{-3}$.
Suppose $c(0)=c_0>0$, then $K=c_0$. Thus,
\begin{equation*}
c(t) = \dfrac{c_0}{(1-k_0 t)^3}.
\end{equation*}
Therefore the free surface is described by the equation $f(t,\sqrt{x^2+y^2})=0$, with
\begin{equation*}
0= f\big(t,\sqrt{x^2+y^2}\big) \vcentcolon= x^2+y^2 - \dfrac{8 c_0}{3 k_0^2 (1-k_0 t)}, \quad \text{where } \ k_0\vcentcolon= \dfrac{\alpha_0^{\prime}}{\gamma\alpha_0}.
\end{equation*}
\paragraph{Time Asymptotic.} The behaviour of free surface $f(t,x^2+y^2)=0$ can be categorised into two cases: $k_0>0$ and $k_0<0$. If $k_0>0$ then the radius of free surface blows up to infinity in finite time as $t\to \frac{1}{k_0}^{-}$. If $k_0<0$, then the radius is shrinking to $0$ as $t\to \infty$.


\section{Flow of Ellipsoid, Hyperboloid, and Cone}\label{sec:john}
We solve (\ref{PDEs}) assuming that potential function takes the following form:
\begin{equation}\label{John}
	\phi(t,x,y,z) = \dfrac{x^2 + y^2 - 2 z^2}{A(t)}, \ \ \text{for some function } \ A(t). 
\end{equation}
This ansatz was first proposed and studied by F. John in \cite{John}.
\subsection{ODE for \texorpdfstring{$A(t)$}{A(t)}}\label{subsec:AODE}
\paragraph{Inferences from the surface evolution equation (\ref{lagrange}).} Suppose that the equation $F(t,\x)=F(t,x,y,z)=0$ describes the moving free boundary surface. For $\x\in\d\Omega(t)$ at time $t\ge0$, let $\{\X(s;t,\x)\}_{0\le s\le t}$ be the backward characteristic curve emanating from $(t,\x)$. Then $F(s,\X(s;t,\x))=0$ for all $0 \le s \le t$. Putting the ansatz (\ref{John}) into the differential equation (\ref{backChar}), one has
\begin{equation*}
	\dfrac{\d X^x}{\d s} = \dfrac{2 X^x }{A(s)}, \qquad \dfrac{\d X^y}{\d s} = \dfrac{2 X^y }{A(s)}, \qquad \dfrac{\d X^z}{\d s} = - \dfrac{4 X^z}{A(s)}.
\end{equation*}
For simplicity, denote $X^{\sigma}_0\equiv X^{\sigma}(0;t,\x)$ for $\sigma=x,\,y,\,z$. Then solving this, one has
\begin{equation*}
	x e^{-g(t)} = X^x_0, \quad y e^{-g(t)} = X^y_0, \quad z e^{2g(t)} = X^z_0, \quad \text{where } \ g(t)\vcentcolon= \int_{0}^{t}\!\!\frac{2}{A(s)}\, \dif s.
\end{equation*}
Suppose that at $s=0$, the initial configuration for the free boundary is determined by a given function $\x \to f(\x) \in \R$, i.e. $F(0,\x)=f(\x)=0$, Then it follows that 
\begin{equation}\label{gS-J}
	0=F(t,\x)= f\big(\X(0;t,\x)\big) = f\big(x e^{-g(t)}, y e^{-g(t)}, z e^{2g(t)}\big).
\end{equation}
Thus the general equations for the free surface takes the form (\ref{gS-J}), with $\x\to f(\x) $ being function independent of the variables $(t,x,y,z)$.

\paragraph{Inferences from the Bernoulli's principle (\ref{bernoulli}).}  Substituting the ansatz (\ref{John}) into the Bernoulli's principle (\ref{bernoulli}), one has
\begin{equation*}
	\Big\{ \dfrac{2-A^{\prime}}{A^2} (x^2+y^2) + \dfrac{2(A^{\prime}+4)}{A^2} z^2 \Big\} \Big\vert_{F(t,x,y,z)=0} = c(t).
\end{equation*}
Rewriting the above equally in accordance to (\ref{gS-J}), we have 
\begin{equation}\label{gB-J}
	\Big\{ \dfrac{2-A^{\prime}}{A^2} e^{2g} \big(|xe^{-g}|^2+|ye^{-g}|^2\big) + \dfrac{2(A^{\prime}+4)}{A^2} e^{-4g} |ze^{2g}|^2 \Big\} \Big\vert_{F(t,x,y,z)=0} = c(t).
\end{equation}
Two equations (\ref{gS-J}) and (\ref{gB-J}) must be consistent, which means that we must have
\begin{equation}\label{consis}
	\dfrac{2-A^{\prime}}{A^2} e^{2g(t)} = c_1 \dfrac{2(A^{\prime}+4)}{A^2}e^{-4g(t)}, \qquad \dfrac{2-A^{\prime}}{A^2} e^{2g(t)} = c_2 c(t),
\end{equation}
for some fixed constants $c_1,\,c_2 \in \R$. The first equation of the above can be rewritten as $6g(t)=\ln(2c_1\frac{4+A^{\prime}}{2-A^{\prime}})$. Taking derivative, and noting that $g(t)=\int_{0}^{t}\frac{2}{A}\dif s$, one gets
\begin{equation}\label{A''}
	\dfrac{A^{\prime\prime}}{4+A^{\prime}}+\dfrac{A^{\prime\prime}}{2-A^{\prime}}=\dfrac{12}{A} \iff A^{\prime\prime} = \dfrac{2(4+A^{\prime})(2-A^{\prime})}{A}.
\end{equation}
Multiplying the above with $A^{\prime}$ we have that $\frac{ A^{\prime\prime}A^{\prime}}{(4+A^{\prime})(2-A^{\prime})}=\frac{2 A^{\prime}}{A}$. Applying the partial fraction decomposition on the left hand side of this equation, we have
\begin{equation*}
	\dfrac{\dif}{\dif t}\ln(4+A^{\prime})^{-\frac{2}{3}}(2-A^{\prime})^{-\frac{1}{3}} =\Big\{\dfrac{-2/3}{4+A^{\prime}} + \dfrac{1/3}{2-A^{\prime}} \Big\} A^{\prime\prime} = \dfrac{2A^{\prime}}{A} = \dfrac{\dif}{\dif t} \ln A^2
\end{equation*}
Multiplying the above with $3$, then integrating the resultant equation we obtain that
\begin{equation}\label{A}
	A^6(2-A^{\prime})(4+A^{\prime})^2 = \mu, \qquad \text{where } \ \mu \vcentcolon= A_0^6 (2-A^{\prime}_0) (4+A^{\prime}_0)^2, 
\end{equation}
with $(A_0,A^{\prime}_0)\equiv(A(0),A^{\prime}(0))$. Next, substituting $6g(t)=\ln(2c_1\frac{4+A^{\prime}}{2-A^{\prime}})$ into the second equation of (\ref{consis}), then using (\ref{A}), we also obtain that for $c_0\equiv c(0)$,
\begin{equation}\label{ct}
	c(t) = \nu (2-A^{\prime})(4+A^{\prime}), \quad \text{where } \ \nu \vcentcolon= \sqrt[\leftroot{-1}\uproot{2}\scriptstyle3]{\dfrac{2c_1}{c_2^3 \mu}} = \dfrac{c_0}{(2-A_0^{\prime})(4+A_0^{\prime})}.
\end{equation}
\subsection{Solution for \texorpdfstring{$A(t)$}{A(t)} and convergence behaviour}
Expanding the identity (\ref{A}), one obtains $(A^{\prime})^3 +6  (A^{\prime})^2 + \frac{\mu}{A^6} - 32 = 0$. Setting $x=A^{\prime}+2$. Then we obtain the depressed cubic equation of the form:
\begin{equation}\label{cubic}
	x^3 + p x + q = 0\qquad \text{where } \ x=A^{\prime}+2, \ \ p=-12,\ \ q=\dfrac{\mu}{A^6}-16.
\end{equation}
We consider the quantity $$D\vcentcolon= \frac{q^2}{4} + \frac{p^3}{27} = \dfrac{\mu}{4A^6}\big(\dfrac{\mu}{A^6}-32\big).$$ By the formula for depressed cubic equations, if $D>0$ then (\ref{cubic}) has 1 real root and 2 complex roots, and the real root in this case can be obtained using Cardano's formula
\begin{equation}\label{cardano}
	\text{if $D>0$,} \quad
	x=\sqrt[\leftroot{-1}\uproot{2}\scriptstyle 3]{-\dfrac{q}{2}+\sqrt{\frac{q^2}{4} + \frac{p^3}{27}}}+\sqrt[\leftroot{-1}\uproot{2}\scriptstyle 3]{-\dfrac{q}{2}-\sqrt{\frac{q^2}{4} + \frac{p^3}{27}}}
\end{equation}
Furthermore, this real root can also be represented in terms of hyperbolic function:
\begin{equation}\label{hyperbolic}
	\text{if $D>0$,} \quad
	x=-2\, \textrm{sgn}(q) \sqrt{\dfrac{-p}{3}} \cosh \Big\{ \dfrac{1}{3}\textrm{arccosh}\Big(\dfrac{-3|q|}{2p}\sqrt{\dfrac{-3}{p}}\Big)  \Big\},
\end{equation}
where $\textrm{sgn}(q)=\frac{|q|}{q}$ denotes the sign function. On the other hand, if $D<0$ then it has 3 distinct real roots, and if $D=0$ then it has 2 distinct real roots with one of them having multiplicity 2. In these cases, the solution is given by Vi\`ete's trigonometric formula
\begin{equation}\label{viete}
	\text{if $D\le0$,} \quad
	x= 2 \sqrt{-\dfrac{p}{3}}\cos\Big\{ \dfrac{1}{3}\arccos\Big( \dfrac{3q}{2p}\sqrt{\dfrac{-3}{p}} \Big) - \dfrac{2\pi k}{3} \Big\}, \quad \text{for } \ k=0,\,1,\,2.
\end{equation}
Setting $a(t)\vcentcolon= A(t)/A_0$. By the definition of $\mu$ in (\ref{A}), it follows that
\begin{equation}\label{Daf}
	D>0 \ \text{ if and only if } \ \dfrac{1}{|a(t)|^6} \vcentcolon= \Big|\dfrac{A_0}{A(t)}\Big|^6 > \dfrac{32}{(2-A_0^{\prime})(4+A_0^{\prime})^2}=\vcentcolon \beta(A_0^{\prime}).
\end{equation}
Suppose that $A(t)$ is continuous near $t=0$. Then for small time $0<t<\!\!< 1$, one has $|a(t)-1|<\!\!<1$. Therefore according to (\ref{Daf}), the root of (\ref{cubic}) on a neighbourhood of $t=0$ is given by Cardano's formula (\ref{cardano}) if $1>\beta(A_0^{\prime})$, and Vi\`ete's formula (\ref{viete}) if $1\le \beta(A_0^{\prime})$. It can be verified that $1 > \beta(y)$ for $y\in(\infty,-6)\cup(2,\infty)$ and $1\le \beta(y)$ for $y\in[-6,-4)\cup(-4,2)$. Therefore we split the analysis into 2 cases:
\paragraph{Case 1: $A_0^{\prime}\in[-6,-4)\cup(-4,2)$.}\label{par:Ncase1} The root for (\ref{cubic}) is given by one of the three possible solutions of Vi\`ete's trigonometric formula (\ref{viete}):
\begin{equation}\label{Vk}
	A^{\prime}(t)+2\!=\! V_{k}(a(t))\!\vcentcolon=\! 4 \cos\Big\{ \dfrac{1}{3}\arccos\Big( 1- \dfrac{2}{|a(t)|^6 \beta(A_0^{\prime})}  \Big) -\dfrac{2\pi k}{3} \Big\},
\end{equation} 
for $k=0,\,1,\,2$. The relevant integer $k$ must be chosen so that solution satisfies the initial condition $ A_0^{\prime}+2 = V_k(1) $. By the definition of $\beta(\cdot)$ in (\ref{Daf}), it can be verified that $\arccos(1-\tfrac{2}{\beta(A_0^{\prime})}) \in [0,\pi]$ for $A_0^{\prime}\in[-6,2]$. Using this, we have
\begin{equation*}
V_2(1)-2\in [-6,-4],\quad V_1(1)-2\in[-4,0],\quad V_0(1)-2\in[0,2], \qquad \text{ for } \ A_0^{\prime}\in[-6,2].
\end{equation*}
Therefore, depending on the initial data $A_0^{\prime}$, $a(t)$ solves the differential equation:
\begin{equation}\label{dadt}
	A_0\dfrac{\dif a}{\dif t} = A^{\prime}(t) = \left\{ \begin{aligned}
		&V_2\big(a(t)\big) - 2 && \text{if } \ A_0^{\prime}\in [-6,-4),\\
		&V_1\big(a(t)\big) - 2 && \text{if } \ A_0^{\prime}\in (-4,0],\\
		&V_0\big(a(t)\big) - 2  && \text{if } \ A_0^{\prime}\in [0,2).
	\end{aligned}\right.
\end{equation}
Integrating the above in $s\in[0,t]$, and by the fact that $a(0)=1$, we have for $t\ge0$, 
\begin{equation}\label{AIntE}
	\mathcal{I}_k\big(a(t)\big) \vcentcolon= \int_{1}^{a(t)} \dfrac{\dif s}{V_k(s)-2} = \dfrac{t}{A_0}, \qquad \text{for } \ \left\{\begin{aligned}
		&k=2 \  \text{ if } \ A_0^{\prime}\in[-6,-4),\\
		&k=1 \  \text{ if } \ A_0^{\prime}\in(-4,0],\\
		&k=0 \  \text{ if } \ A_0^{\prime}\in[0,2).
	\end{aligned}\right.
\end{equation}
Let $T>0$ be the maximal time for which $a(t)\mapsto \mathcal{I}_k(a(t))$ is invertible. Then
\begin{equation}\label{AIntE-1}
	A(t) = A_0 \mathcal{I}_k^{-1} \big( \dfrac{t}{A_0} \big) \qquad \text{for } \ k=0,\,1,\,2 \ \text{ and } \ t\in[0,T).
\end{equation}
\paragraph{Case 2: $A_0^{\prime}\in(-\infty,-6]\cup(2,\infty)$.}\label{par:Ncase2} In this case, the one physical real root for (\ref{cubic}) is given by the Cardano's formula (\ref{cardano}):
\begin{subequations}\label{W}
	\begin{align}
		&A^{\prime}(t) = W(a(t)) \vcentcolon= W_{+}\big(a(t)\big)+W_{-}\big(a(t)\big),\label{Wa}\\
		\text{where } \  &W_{\pm}(a) \vcentcolon= 2 \Big(1-\dfrac{2}{a^6 \beta(A_0^{\prime})} \big\{ 1 \pm\sqrt{1-a^6 \beta(A_0^{\prime})} \big\} \Big)^{\frac{1}{3}}.\label{Wb}
	\end{align}
	We remark that $W$ can also be represented in terms of hyperbolic function
	\begin{equation}\label{Wc}
		W(a) = -4 \, \textrm{sgn}\Big(\dfrac{2}{a^6 \beta(A_0^{\prime})} - 1 \Big) \cosh \Big\{ \dfrac{1}{3} \textrm{arccosh} \Big( \big| 1- \dfrac{2}{|a(t)|^6 \beta(A_0^{\prime})} \big| \Big) \Big\}.
	\end{equation}
\end{subequations}
Let $T\vcentcolon=\sup\{ t>0 \vcentcolon |a(t)|^6 > \beta(A_0^{\prime}) \}$, then $a(t)$ satisfies the differential equation
\begin{equation}\label{dadtW}
	A_0 \dfrac{\dif a}{\dif t} = A^{\prime}(t) = W\big( a(t) \big)-2, \qquad \text{for } \ t\in[0,T).
\end{equation}
Integrating the above in $s\in[0,t]$, and using $a(0)=1$ we have that for $t\in[0,T)$,
\begin{equation}\label{AIntH}
	\mathcal{I}_H\big(a(t)\big) \vcentcolon= \int_{1}^{a(t)} \dfrac{\dif s}{W(s)-2} = \dfrac{t}{A_0} \ \Rightarrow \ A(t) = A_0 \mathcal{I}_H^{-1}\big( \dfrac{t}{A_0} \big).
\end{equation}
Using the integral expressions (\ref{AIntE}) and (\ref{AIntH}), one can numerically plot the curve $a\mapsto \mathcal{I}_0(a),\mathcal{I}_1(a),\mathcal{I}_2(a),\mathcal{I}_H(a)$, which is then analysed to determine the convergence of $(A(t),A^{\prime}(t))$ as $t\to\pm\infty$. Without loss of generality, we restrict our analysis to the case $A_0>0$. The case $A_0<0$ can be obtained by the transformation: $(t,A)\mapsto (-t,-A)$, which preserves the differential equations (\ref{A''})--(\ref{A}). We consider initial data $(A_0,A_0^{\prime})$ in $5$ distinct regimes. The statements and figures listed below, describing the behaviour of $(A,A^{\prime})$ as $t\to\infty$ are based on the numerical integration of (\ref{AIntE}) and (\ref{AIntH}).

\paragraph{Regime 1: $A_0^{\prime}\in ( -\infty,-6]$.}\label{par:regime1} In this case, $1\ge \beta(A_0^{\prime})$ and the differential equation is given by $A^{\prime}(t)=W(\frac{A(t)}{A_0})-2$. Moreover, (\ref{A}) implies that $\mu>0$ and $t\mapsto A^{\prime}(t)$ is monotone decreasing due to $A>0$ and (\ref{A''}). From the fact that $\mu>0$ and (\ref{A}), it follows that as $t$ increases from $0$, one has the limit $A\to 0^{+}$ and $A^{\prime}\to-\infty$. On the other hand, as $t$ decreases from $0$, one has the limit $A\to A_{\ast}^{-}$ and $A^{\prime}\to -6^{-}$ where $A_{\ast}\vcentcolon=(\frac{\mu}{(2-A^{\prime})(4+A^{\prime})^2})^{\frac{1}{6}}\vert_{A^{\prime}=-6}$. Thus, we aim to find the times: $T_{+}>0$ for which $A^{\prime}\to -\infty$, and $T_{-}<0$ for which $A^{\prime}\to -6$. By the numerical simulation, we have
\begin{enumerate}[label=(\roman*),ref=(\roman*),font=\textnormal]
	\item\label{item:reg1-i} For the case of increasing $t>0$, the integral (\ref{AIntH}) converges to a positive finite time: $\mathcal{I}_{H}(a)\to \frac{T_{+}}{A_0}>0$ as $a=\frac{A}{A_0}\to 0$. Thus $(A,A^{\prime})\to (0,-\infty)$ in finite time as $t\to T_{+}$.
	\item\label{item:reg1-ii} For the case of decreasing $t<0$, the integral (\ref{AIntH}) converges to a negative finite time: $\mathcal{I}_{H}(a)\to \frac{T_{-}}{A_0}<0$ as $a=\frac{A}{A_0}\to \frac{A_{\ast}}{A_0}$. Thus $(A,A^{\prime})\to (A_{\ast},-6)$ in finite time as $t\to T_{-}$. In this scenario, we take $(A,A^{\prime})\vert_{t=T_{-}}=(A_{\ast},-6)$ as an initial data starting in \hyperref[item:reg2-ii]{\textbf{Regime 2(ii)}} then concatenate the resulting solution with $A(t)\vert_{ T_{-} \le t\le 0}$.
\end{enumerate}
\begin{figure}[H]
	\centering
	\includegraphics[trim={1.15cm 0 0 0},clip,width=\textwidth]{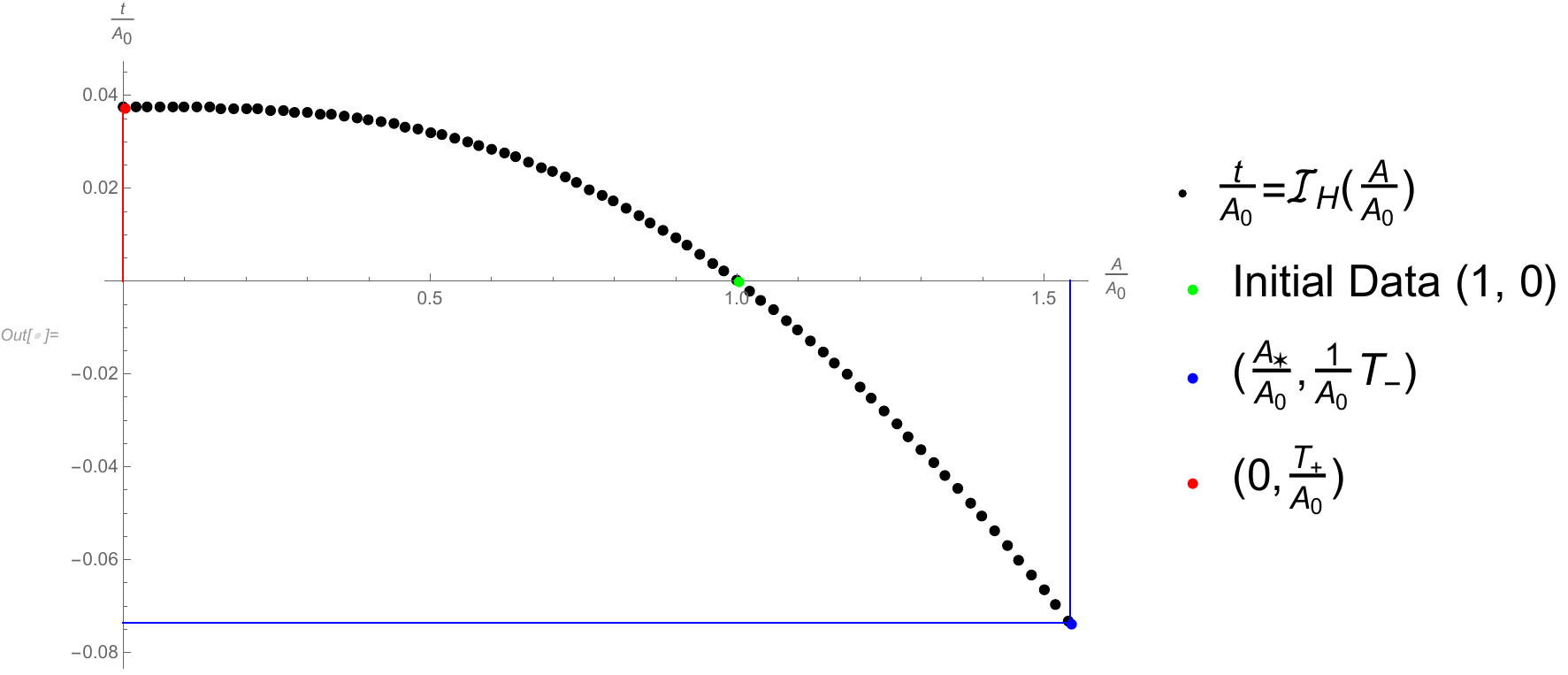}
	\caption{the plot of curve $(\frac{t}{A_0},\mathcal{I}_H(\frac{A}{A_0}))$ with $A_0^{\prime}=-10$ for Regime 1.}
	\label{fig:Regime1}
\end{figure}

\paragraph{Regime 2: $A_0^{\prime}\in [-6,-4)$.}\label{par:regime2} In this case, $1< \beta(A_0^{\prime})$ and $A$ solves $A^{\prime}(t)=V_2(\frac{A(t)}{A_0})-2$ according to (\ref{dadt}). Moreover, (\ref{A}) implies that $\mu>0$, and $t\mapsto A^{\prime}(t)$ is monotone decreasing since $A>0$ and (\ref{A''}). From the fact that $\mu>0$ and (\ref{A}), it follows that as $t$ increases from $0$, one has the limit $A\to A_{\ast}^{+}$ and $A^{\prime}\to-6^{+}$, where $A_{\ast}\vcentcolon=(\frac{\mu}{(2-A^{\prime})(4+A^{\prime})^2})^{\frac{1}{6}}\vert_{A^{\prime}=-6}$. On the other hand, as $t$ decreases from $0$, one has the limit $A\to \infty$ and $A^{\prime}\to -4^{-}$. Thus, we aim to find the times $T_{+}>0$ for which $A^{\prime}\to -6$, and $T_{-}<0$ for which $A^{\prime}\to -4^{-}$.  According to the numerical simulation, we have
\begin{enumerate}[label=(\roman*),ref=(\roman*),font=\textnormal]
	\item\label{item:reg2-i} For the case of increasing $t>0$, the integral (\ref{AIntE}) converges to a positive finite time: $\mathcal{I}_{2}(a)\to \frac{T_{+}}{A_0}>0$ as $a=\frac{A}{A_0}\to \frac{A_{\ast}}{A_0}$. Thus $(A,A^{\prime})\to (A_{\ast},-6)$ in finite time as $t\to T_{+}$. In this scenario, we take $(A,A^{\prime})\vert_{t=T_{+}}=(A_{\ast},-6)$ as an initial data starting in \hyperref[item:reg1-i]{\textbf{Regime 1(i)}} then concatenate the resulting solution with $A(t)\vert_{ 0 \le t\le T_{+}}$.
	\item\label{item:reg2-ii} For the case of decreasing $t<0$, it can be verified that $\frac{1}{V_2(a)-2}\to -\frac{1}{4}$ as $a\to\infty$. Hence the integral (\ref{AIntE}) diverges to negative infinity: $\mathcal{I}_{2}(a)\to -\infty$ as $a=\frac{A}{A_0}\to \infty$, which implies that $T_{-}=-\infty$. Thus $(A,A^{\prime})\to (\infty,-4)$ in infinite time as $t\to -\infty$. 
\end{enumerate}
\begin{figure}[H]
	\centering
	\includegraphics[trim={1.15cm 0 0 0},clip,width=\textwidth]{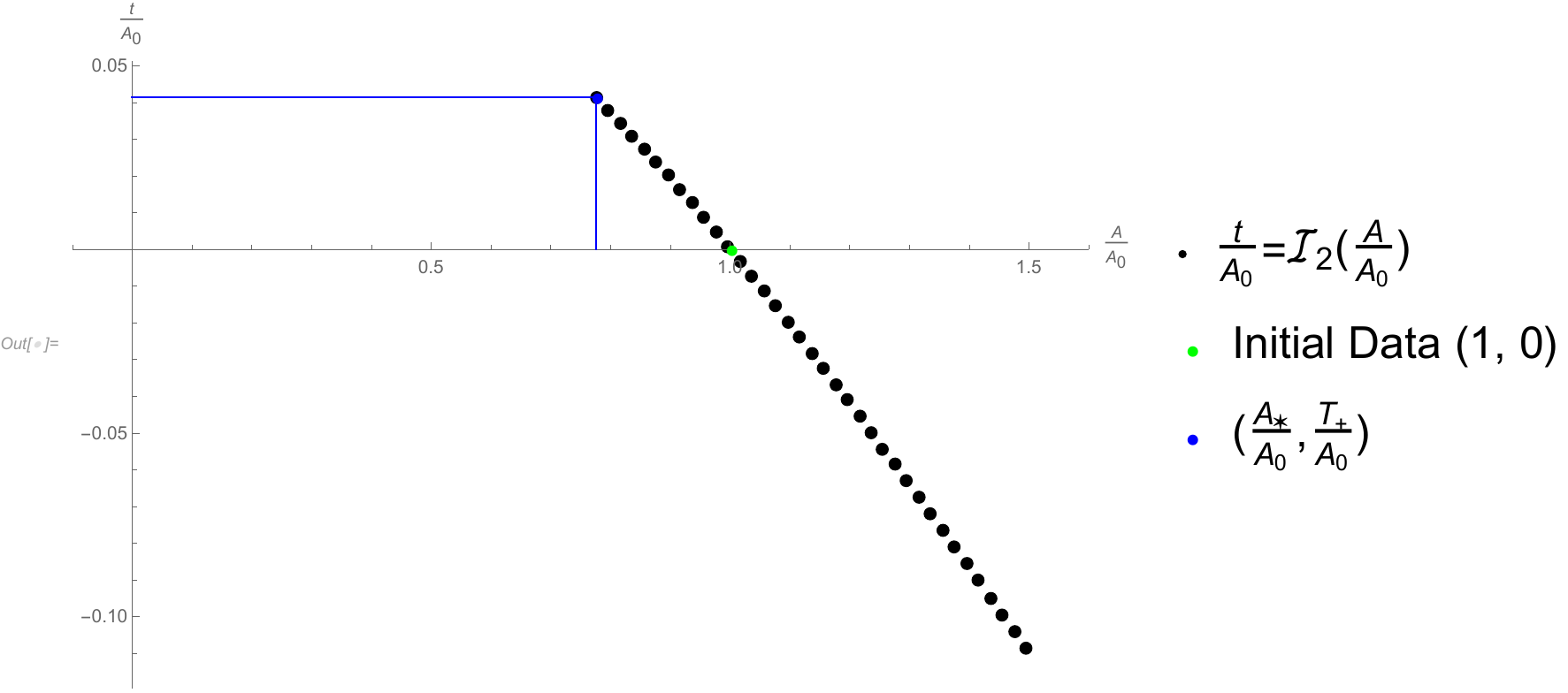}
	\caption{the plot of curve $(\frac{t}{A_0},\mathcal{I}_2(\frac{A}{A_0}))$ with $A_0^{\prime}=-5$ for Regime 2.}
	\label{fig:Regime2}
\end{figure}
\begin{figure}[H]
	\centering
	\includegraphics[trim={1.15cm 0 0 0},clip,width=\textwidth]{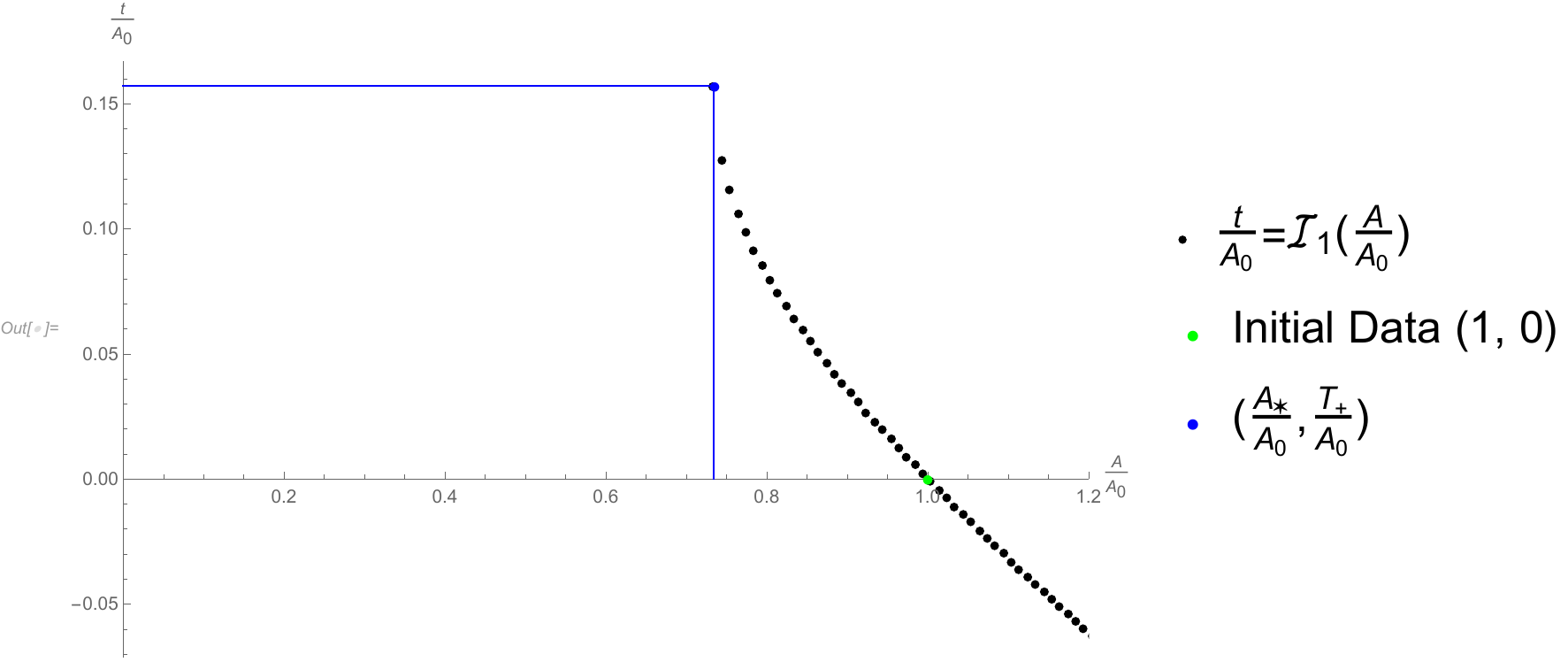}
	\caption{the plot of curve $(\frac{t}{A_0},\mathcal{I}_1(\frac{A}{A_0}))$ with $A_0^{\prime}=-3$ for Regime 3.}
	\label{fig:Regime3}
\end{figure}

\paragraph{Regime 3: $A_0^{\prime}\in (-4,0]$.}\label{par:regime3} In this case, $1< \beta(A_0^{\prime})$ and the differential equation is given by $A^{\prime}(t)=V_1(\frac{A(t)}{A_0})-2$ according to (\ref{dadt}). Moreover, (\ref{A}) implies that $\mu>0$, and $t\mapsto A^{\prime}(t)$ is monotone increasing since $A>0$ and (\ref{A''}). From the fact that $\mu>0$ and (\ref{A}), it follows that as $t$ increases from $0$, one has the limit $A\to A_{\ast}^{+}$ and $A^{\prime}\to 0^{-} $, where $A_{\ast}\vcentcolon=(\frac{\mu}{(2-A^{\prime})(4+A^{\prime})^2})^{\frac{1}{6}}\vert_{A^{\prime}=0}$. On the other hand, as $t$ decreases from $0$, one has the limit $A\to \infty$ and $A^{\prime}\to -4^{+}$. Thus, we aim to find the times $T_{+}>0$ for which $A^{\prime}\to 0$, and $T_{-}<0$ for which $A^{\prime}\to-4$. According to the numerical simulation, we have 
\begin{enumerate}[label=(\roman*),ref=(\roman*),font=\textnormal]
	\item\label{item:reg3-i} For the case of increasing $t>0$, the integral (\ref{AIntE}) converges to a positive finite time: $\mathcal{I}_{1}(a)\to \frac{T_{+}}{A_0}>0$ as $a=\frac{A}{A_0}\to \frac{A_{\ast}}{A_0}$. Thus $(A,A^{\prime})\to (A_{\ast},0)$ in finite time as $t\to T_{+}$. In this scenario, we take $(A,A^{\prime})\vert_{t=T_{+}}=(A_{\ast},0)$ as an initial data starting in \hyperref[item:reg4-i]{\textbf{Regime 4(i)}} then concatenate the resulting solution with $A(t)\vert_{ 0 \le t\le T_{+}}$.
	\item\label{item:reg3-ii} For the case of decreasing $t<0$, it can be verified that $\frac{1}{V_1(a)-2}\to -\frac{1}{4}$ as $a\to\infty$. Hence the integral (\ref{AIntE}) diverges to negative infinity: $\mathcal{I}_{1}(a)\to -\infty$ as $a=\frac{A}{A_0}\to \infty$, which implies that $T_{-}=-\infty$. Thus $(A,A^{\prime})\to (\infty,-4)$ as $t\to -\infty$. 
\end{enumerate}

\paragraph{Regime 4: $A_0^{\prime}\in [0,2)$.}\label{par:regime4} In this case, $1< \beta(A_0^{\prime})$ and the differential equation is given by $A^{\prime}(t)=V_0(\frac{A(t)}{A_0})-2$ according to (\ref{dadt}). Moreover, (\ref{A}) implies that $\mu>0$, and $t\mapsto A^{\prime}(t)$ is monotone increasing since $A>0$ and (\ref{A''}). From the fact that $\mu>0$ and (\ref{A}), it follows that as $t$ increases from $0$, one has the limit $A\to \infty$ and $A^{\prime}\to 2^{-} $. On the other hand, as $t$ decreases from $0$, one has the limit $A\to A_{\ast}^{+}$ and $A^{\prime}\to 0^{+}$, where $A_{\ast}\vcentcolon=(\frac{\mu}{(2-A^{\prime})(4+A^{\prime})^2})^{\frac{1}{6}}\vert_{A^{\prime}=0}$. Thus, we aim to find the times $T_{+}>0$ for which $A^{\prime}\to 2^{-}$, and $T_{-}<0$ for which $A^{\prime}\to 0^+$. By the numerical simulation, we have
\begin{enumerate}[label=(\roman*),ref=(\roman*),font=\textnormal]
	\item\label{item:reg4-i} For the case of increasing $t>0$, it can be verified that $\frac{1}{V_0(a)-2}\to \frac{1}{2}$ as $a\to\infty$. Hence the integral (\ref{AIntE}) diverges to positive infinity: $\mathcal{I}_{0}(a)\to \infty$ as $a=\frac{A}{A_0}\to \infty$, which implies that $T_{+}=\infty$. Thus $(A,A^{\prime})\to (\infty,2)$ in infinite time as $t\to \infty$.
	\item\label{item:reg4-ii} For the case of decreasing $t<0$, the integral (\ref{AIntE}) converges to a negative finite time: $\mathcal{I}_{0}(a)\to \frac{T_{-}}{A_0}<0$ as $a=\frac{A}{A_0}\to \frac{A_{\ast}}{A_0}$. Thus $(A,A^{\prime})\to (A_{\ast},0)$ in finite time as $t\to T_{-}$. In this scenario, we take $(A,A^{\prime})\vert_{t=T_{-}}=(A_{\ast},0)$ as an initial data starting in \hyperref[item:reg3-ii]{\textbf{Regime 3(ii)}} then concatenate the resulting solution with $A(t)\vert_{T_{-}\le t\le 0}$.
\end{enumerate}
\begin{figure}[H]
	\centering
	\includegraphics[trim={1.15cm 0 0 0},clip,width=\textwidth]{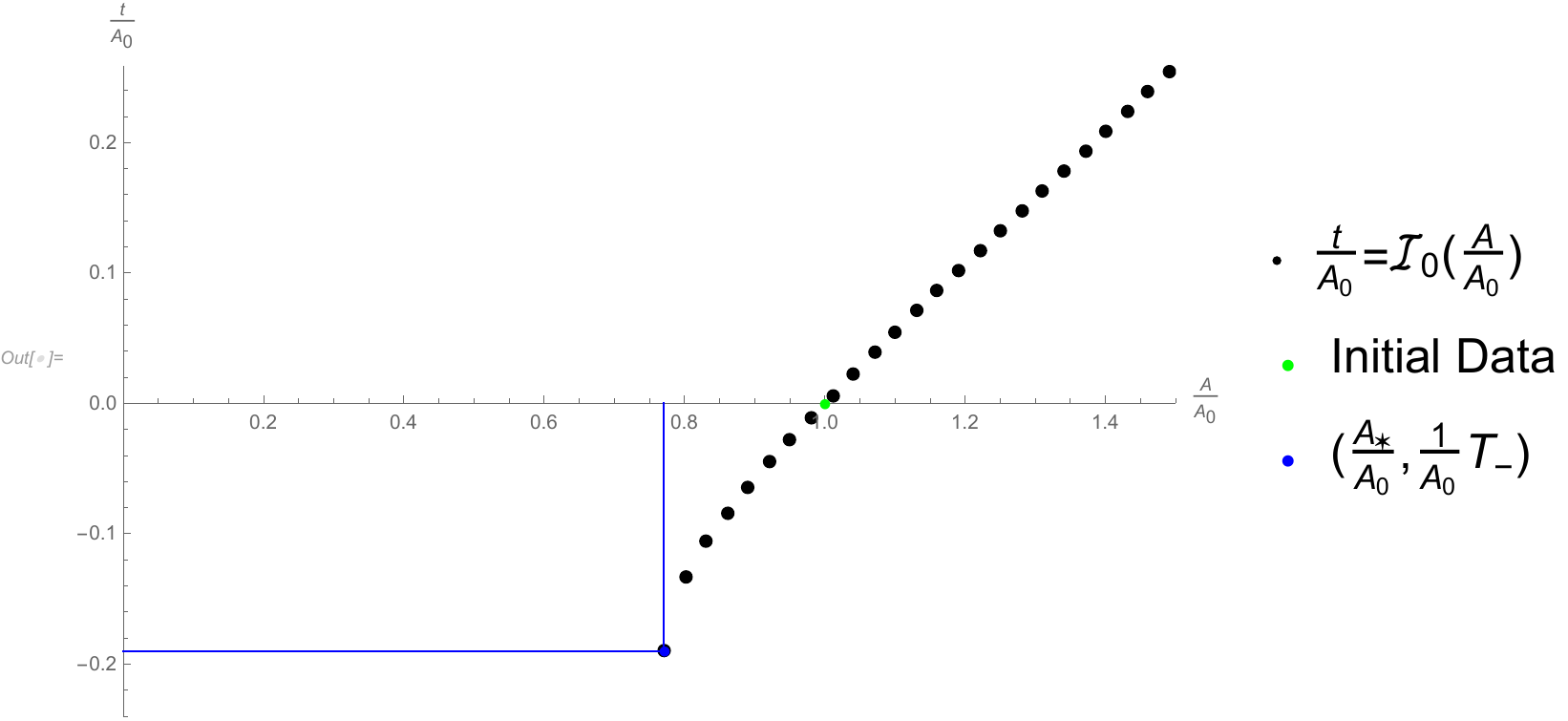}
	\caption{the plot of curve $(\frac{t}{A_0},\mathcal{I}_0(\frac{A}{A_0}))$ with $A^{\prime}_0=1.8$ for Regime 4.}
	\label{fig:Regime4}
\end{figure}

\paragraph{Regime 5: $A_0^{\prime}\in (2,\infty)$.}\label{par:regime5} In this case, $1> \beta(A_0^{\prime})$ and the differential equation is given by $A^{\prime}(t)=W(\frac{A(t)}{A_0})-2$. Moreover, (\ref{A}) implies that $\mu<0$, and $t\mapsto A^{\prime}(t)$ is monotone decreasing since $A>0$ and (\ref{A''}). From the fact that $\mu<0$ and (\ref{A}), it follows that as $t$ increases from $0$, one has the limit $A\to \infty$ and $A^{\prime}\to 2^{+} $. On the other hand, as $t$ decreases from $0$, one has the limit $A\to 0^+$ and $A^{\prime}\to \infty$. Thus, we aim to find the time $T_{+}$ for which $A^{\prime}\to 2^{+}$. From (\ref{A}), we see that $A\to\infty$ as $A^{\prime}\to 2^{+}$, and $T_{-}<0$ for which $A^{\prime}\to \infty$. By the numerical simulation, we have
\begin{enumerate}[label=(\roman*),ref=(\roman*),font=\textnormal]
	\item\label{item:reg5-i} For the case of increasing $t>0$, it can be verified that $\frac{1}{W(a)-2}\to \frac{1}{2}$ as $a\to\infty$. Hence the integral (\ref{AIntH}) diverges to positive infinity: $\mathcal{I}_{H}(a)\to \infty$ as $a=\frac{A}{A_0}\to \infty$, which implies that $T_{+}=\infty$. Thus $(A,A^{\prime})\to (\infty,2)$ in infinite time as $t\to \infty$.
	\item\label{item:reg5-ii} For the case of decreasing $t<0$, the integral (\ref{AIntH}) converges to a negative finite time: $\mathcal{I}_{0}(a)\to \frac{T_{-}}{A_0}<0$ as $a=\frac{A}{A_0}\to 0$. Thus $(A,A^{\prime})\to (0,\infty)$ in finite time as $t\to T_{-}$.
\end{enumerate}
\begin{figure}[H]
	\centering
	\includegraphics[trim={1.15cm 0 0 0},clip,width=\textwidth]{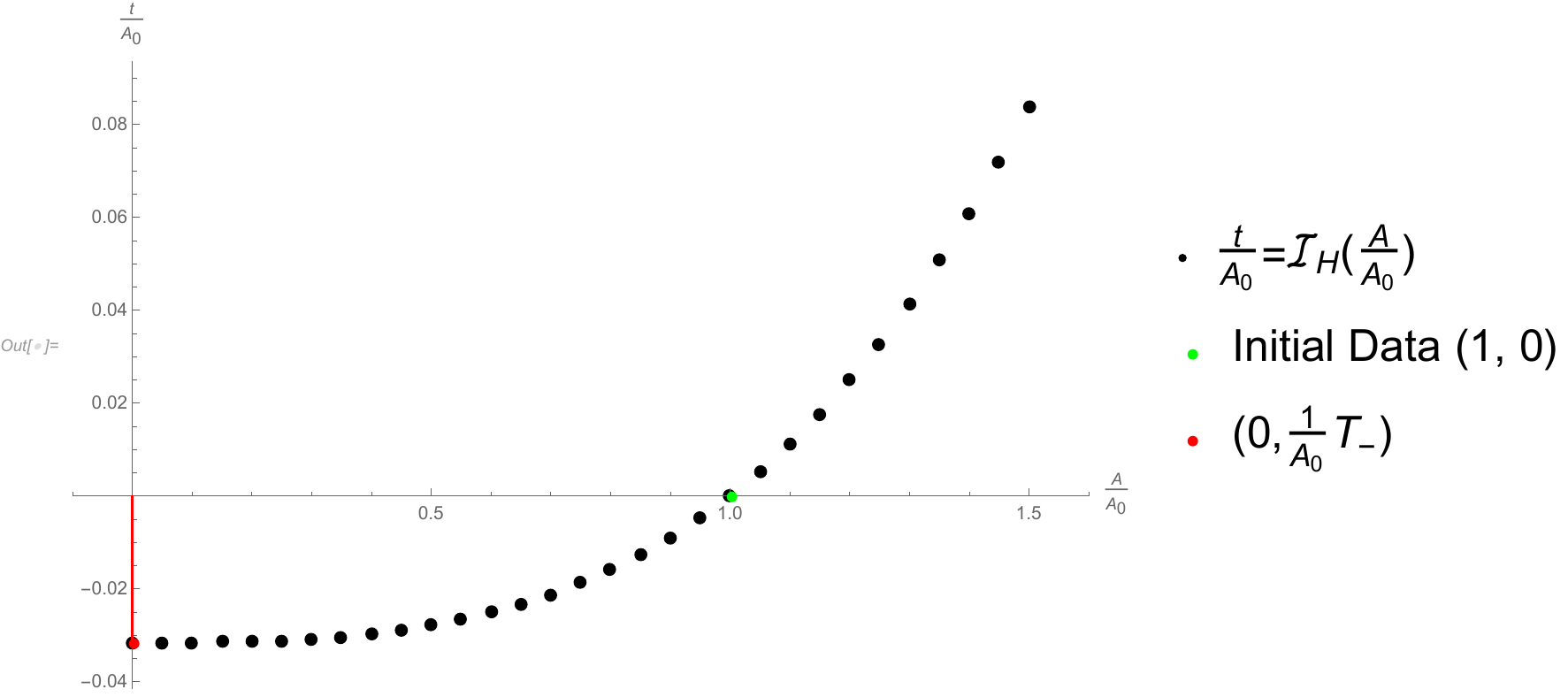}
	\caption{the plot of curve $(\frac{t}{A_0},\mathcal{I}_H(\frac{A}{A_0}))$ with $A_0^{\prime}=10$ for Regime 5.}
	\label{fig:Regime5}
\end{figure}

\subsection{Free boundary surface and phase space}\label{subsec:phase}
The free boundary surface of the flow is determined by (\ref{gB-J}), (\ref{A}), and (\ref{ct}). This is stated as:
\begin{equation}\label{sJohn}
	(2-A^{\prime})(4+A^{\prime}) (x^2+y^2) + 2 (4+A^{\prime})^2 z^2 = \dfrac{\mu\nu}{A^4},
\end{equation}
with $\mu=A_0^6(2-A_0^{\prime})(4+A_0^{\prime})^2$ and $\nu=\frac{c_0}{(2-A_0^{\prime})(4+A_0^{\prime})}$. From this, we see that
\begin{itemize}
	\item If $A^{\prime}< -4$ or $A^{\prime}>2$, then the surface is either an one-sheeted or two-sheeted hyperboloid, or a cone;
	\item If $-4<A^{\prime}<2$, then the surface is an ellipsoid.
\end{itemize}
For a given initial data $A_0=A(0)$ and $A_0^{\prime}=A^{\prime}(0)$, the solution to (\ref{A''})--(\ref{A}) exists up to the maximum time of existence $T\neq 0$ and $A(t)\in \mathcal{C}^2\big([0,T)\big)$ or $A(t)\in \mathcal{C}^2\big((T,0]\big)$ due to the construction given in the previous section. In what follows, we study the evolution of its corresponding free surface. Once again, due to the time symmetry of the differential equations (\ref{A''})--(\ref{A}) for $(t,A) \leftrightarrow (-t,-A)$, we restrict our analysis to the case $A_0>0$. We split the initial data into 9 cases:

\paragraph{Case \textrm{I}(a): $A_0^{\prime}<-4$ and $c_0>0$. The initial free boundary is an One-Sheeted Hyperboloid.}\label{par:case1a} In this case, $\mu>0$ and $\nu<0$. Rewriting (\ref{sJohn}) using (\ref{A}), we have
\begin{equation}\label{hyp}
	x^2+y^2 + \dfrac{2(4+A^{\prime})}{2-A^{\prime}} z^2 = \nu \sqrt[\leftroot{-1}\uproot{2}\scriptstyle 3]{\dfrac{\mu (4+A^{\prime})}{2-A^{\prime}}}.
\end{equation}
By the analysis of \hyperref[par:regime1]{\textbf{Regime 1}} and \hyperref[par:regime2]{\textbf{Regime 2}} in the previous section, we have
\begin{enumerate}[label=(\roman*),font=\textnormal]
	\item For increasing $t>0$, there exists a finite time $T_{+}>0$ such that $A(t)\to 0$, $A^{\prime}(t) \to -\infty$ as $t\to T_{+}$. Thus, as $t\to T_{+}$, the free surface converges to an \textbf{One-Sheeted Hyperboloid} given by the following equation:
	\begin{equation*}
		x^2+y^2 -2 z^2 = -\nu \sqrt[\leftroot{-1}\uproot{2}\scriptstyle 3]{\mu} > 0.
	\end{equation*}
	\item For decreasing $t<0$, we have the limit $A(t)\to \infty$, $A^{\prime}(t) \to -4$ as $t\to -\infty$. Thus, as $t\to -\infty$, the free surface collapses to the \textbf{$z$-axis} described by the equation $x^2+y^2 =0$.
\end{enumerate}

\paragraph{Case \textrm{I}(b): $A_0^{\prime}<-4$ and $c_0<0$. The initial free boundary is a Two-Sheeted Hyperboloid.}\label{par:case1b} In this case, $\mu>0$ and $\nu>0$. The evolution of surface is described by the same equation as (\ref{hyp}). According to \hyperref[par:regime1]{\textbf{Regime 1}} and \hyperref[par:regime2]{\textbf{Regime 2}}, we have
\begin{enumerate}[label=(\roman*),font=\textnormal]
	\item For increasing $t>0$, there exists a finite time $T_{+}>0$ such that $A(t)\to 0$, $A^{\prime}(t) \to -\infty$ as $t\to T_{+}$. Thus, as $t\to T_{+}$, the free surface converges to a \textbf{Two-Sheeted Hyperboloid} given by the following equation
	\begin{equation*}
		x^2+y^2 -2 z^2 = -\nu \sqrt[\leftroot{-1}\uproot{2}\scriptstyle 3]{\mu} < 0.
	\end{equation*}
	\item For decreasing $t<0$, we have the limit $A(t)\to \infty$, $A^{\prime}(t) \to -4$ as $t\to -\infty$. Thus, as $t\to -\infty$, the free surface collapses to the \textbf{$z$-axis} described by the equation $x^2+y^2 =0$.
\end{enumerate}

\paragraph{Case \textrm{I}(c): $A_0^{\prime}<-4$ and $c_0=0$. The initial free boundary is a Cone.}\label{par:case1c} In this case, $\mu>0$ and $\nu=0$. The evolution of surface is described by the same equation as (\ref{hyp}). By the analysis of \hyperref[par:regime1]{\textbf{Regime 1}} and \hyperref[par:regime2]{\textbf{Regime 2}} in the previous section, we have
\begin{enumerate}[label=(\roman*),font=\textnormal]
	\item For increasing $t>0$, there exists a finite time $T_{+}>0$ such that $A(t)\to 0$, $A^{\prime}(t) \to -\infty$ as $t\to T_{+}$. Thus, as $t\to T_{+}$, the free surface converges to a \textbf{Cone} with slope $2$ given by the equation
	\begin{equation*}
		x^2+y^2 -2 z^2 = 0.
	\end{equation*}
	\item For decreasing $t<0$, we have the limit $A(t)\to \infty$, $A^{\prime}(t) \to -4$ as $t\to -\infty$. Thus, as $t\to -\infty$, the free surface collapses to the \textbf{$z$-axis} described by the equation $x^2+y^2 =0$.
\end{enumerate}

\paragraph{Case \textrm{II}: $-4<A_0^{\prime}<2$ and $c_0\ge0$. The initial free boundary is an Ellipsoid or a point.}\label{par:case2} The free surface is described by (\ref{sJohn}) and solution $A(t)$ as
\begin{equation*}
(2-A^{\prime})(4+A^{\prime}) (x^2+y^2) + 2 (4+A^{\prime})^2 z^2 = \dfrac{\mu\nu}{A^4}.
\end{equation*}
If $c_0=0$ then (\ref{sJohn}) reduces to the equation $(2-A^{\prime})(x^2+y^2)+2(4+A^{\prime})z^2 = 0$ which describes a point at the origin. In addition, we remark that $c_0<0$ is not possible in this case. By the analysis of \hyperref[par:regime3]{\textbf{Regime 3}} and \hyperref[par:regime4]{\textbf{Regime 4}} in the previous section, we have
\begin{enumerate}[label=(\roman*),font=\textnormal]
	\item For increasing $t>0$, we have $A(t)\to \infty$, $A^{\prime}(t) \to 2$ as $t\to \infty$. Therefore in the limit $t\to\infty$, the equation (\ref{sJohn}) reduces to $z^2=0$, which means that the free surface collapses into the \textbf{$(x,y)$-Plane}.
	\item For decreasing $t<0$, we have the limit $A(t)\to \infty$, $A^{\prime}(t) \to -4$ as $t\to -\infty$. Thus, as $t\to -\infty$, the free surface collapses to the \textbf{$z$-Axis} described by the equation $x^2+y^2 =0$.
\end{enumerate}

\paragraph{Case \textrm{III}(a): $A_0^{\prime}>2$ and $c_0>0$. The initial free boundary is a Two-Sheeted Hyperboloid.}\label{par:case3a} In this case, $\mu<0$ and $\nu<0$. The free surface is described by the same equation as (\ref{hyp}). By the analysis of \hyperref[par:regime5]{\textbf{Regime 5}} in the previous section, we have
\begin{enumerate}[label=(\roman*),font=\textnormal]
	\item For increasing $t>0$, we have $A(t)\to \infty$, $A^{\prime}(t) \to 2^{+}$ as $t\to \infty$. Taking this limit in (\ref{sJohn}), the equation reduces to $z^2=0$ as $t\to\infty$. Thus, the free surface collapses to the \textbf{(x,y)-Plane} in the limit $t\to \infty$.
	\item For decreasing $t<0$, there exists a finite time $-\infty<T_{-}<0$ such that $A(t)\to 0$, $A^{\prime}(t) \to \infty$ as $t\to T_{-}$. Thus, as $t\to T_{-}$, the free surface converges to a \textbf{Two-Sheeted Hyperboloid} described by the equation:
	\begin{equation*}
		x^2+y^2-2z^2 = -\nu \sqrt[\leftroot{-1}\uproot{2}\scriptstyle 3]{\mu}<0
	\end{equation*}
\end{enumerate}

\paragraph{Case \textrm{III}(b): $A_0^{\prime}>2$ and $c_0<0$. The initial free boundary is an One-Sheeted Hyperboloid.}\label{par:case3b} In this case, $\mu<0$ and $\nu>0$. The free surface is described by the same equation as (\ref{hyp}). According to \hyperref[par:regime5]{\textbf{Regime 5}} in the previous section, we have
\begin{enumerate}[label=(\roman*),font=\textnormal]
	\item For increasing $t>0$, we have $A(t)\to \infty$, $A^{\prime}(t) \to 2^{+}$ as $t\to \infty$. Taking this limit in (\ref{sJohn}), the equation reduces to $z^2=0$ as $t\to\infty$. Thus, the free surface collapses to the \textbf{(x,y)-Plane} in the limit $t\to \infty$.
	\item For decreasing $t<0$, there exists a finite time $-\infty<T_{-}<0$ such that $A(t)\to 0$, $A^{\prime}(t) \to \infty$ as $t\to T_{-}$. Thus, as $t\to T_{-}$, the free surface converges to an \textbf{One-Sheeted Hyperboloid} described by the equation:
	\begin{equation*}
		x^2+y^2-2z^2 = -\nu \sqrt[\leftroot{-1}\uproot{2}\scriptstyle 3]{\mu}>0
	\end{equation*}
\end{enumerate}
\paragraph{Case \textrm{III}(c): $A_0^{\prime}>2$ and $c_0=0$. The initial free boundary is a Cone.}\label{par:case3c} In this case, $\mu<0$ and $\nu=0$. The free surface is an evolving cone described by the same equation as (\ref{hyp}). By the analysis in \hyperref[par:regime5]{\textbf{Regime 5}}, we have
\begin{enumerate}[label=(\roman*),font=\textnormal]
	\item For increasing $t>0$, we have $A(t)\to \infty$, $A^{\prime}(t) \to 2^{+}$ as $t\to \infty$. Taking this limit in (\ref{sJohn}), the equation reduces to $z^2=0$ as $t\to\infty$. Thus, the free surface collapses to the \textbf{(x,y)-Plane} in the limit $t\to \infty$.
	\item For decreasing $t<0$, there exists a finite time $-\infty<T_{-}<0$ such that $A(t)\to 0$, $A^{\prime}(t) \to \infty$ as $t\to T_{-}$. Thus, as $t\to T_{-}$, the free surface converges to the \textbf{Cone} with slope $2$ described by the equation:
	\begin{equation*}
		x^2+y^2-2z^2 = 0
	\end{equation*}
\end{enumerate}
\paragraph{Case \textrm{IV}: $A_0\in\R$ and $A_0^{\prime}=-4$. The initial free boundary is the $z$-axis.}\label{par:case4} For this initial data, $\mu\nu = 0$ and the solution is explicitly given by $A(t)=-4 t + A_0$. The equation (\ref{sJohn}) for the surface remains as the \textbf{$z$-Axis} given by the equation $x^2+y^2=0$.

\paragraph{Case \textrm{V}: $A_0\in\R$ and $A_0^{\prime}=2$. The initial free boundary is a Two-Sheeted Plane perpendicular to the $z$-axis.}\label{par:case5} For this initial data, $\mu=0$ and $\mu\nu = 6 c_0 A_0^6$. The solution is explicitly given by $A(t)=2t+A_0$, and the equation (\ref{sJohn}) for the surface is:
\begin{equation*}
	z^2 = \dfrac{c_0 |A(t)|^6}{9} = \dfrac{c_0}{9} |2t+A_0|^2.
\end{equation*}
This describes a two-sheeted plane perpendicular to the $z$-axis moving along the $z$-direction with velocities $\pm\frac{c_0(2t+A_0)}{9}$. Note that in this case $c_0<0$ is impossible. 
\begin{figure}[H]
	\includegraphics[width=\textwidth,trim={1.15cm 0 0 0},clip]{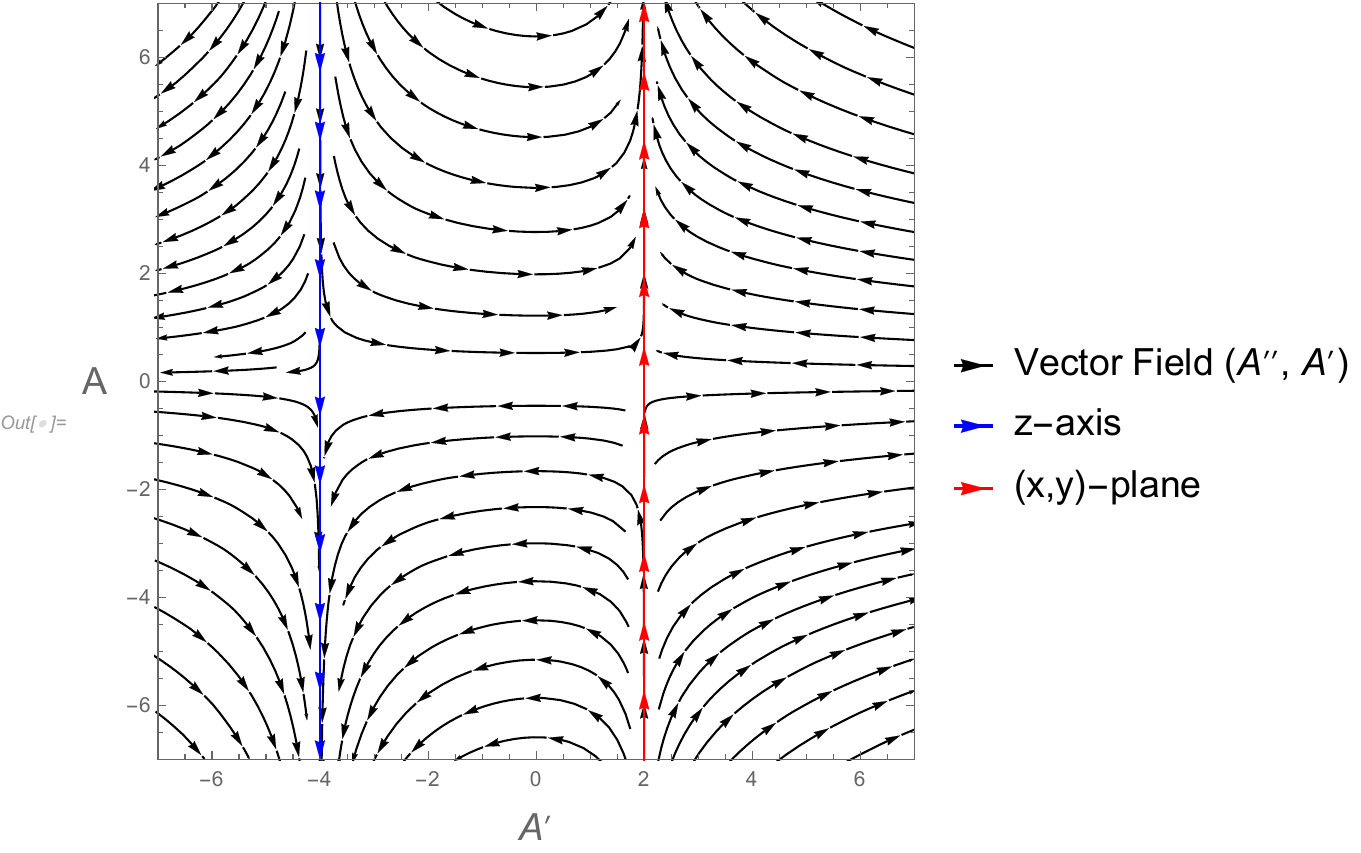}
\end{figure}
\begin{remark}
	The point $(A_0^{\prime},A_0)$ for which $A_0^{\prime}(t)\notin \{-4,2\}$ and $A_0=0$ cannot exist as a solution to (\ref{A}). Thus the set $\{ (q,p)\in\R^2 \vcentcolon q\in(-\infty,-4)\cup(-4,2)\cup(2,\infty) \text{ and } p=0 \}$ is excluded from the phase space of $(A^{\prime},A)$.
\end{remark}

\subsection{Time asymptotic of \texorpdfstring{$(A,A^{\prime})$}{(A,A')}}
In this section, we study the time asymptotic behaviour of $A(t)$ and $A^{\prime}(t)$. As before, we restrict our study to the case $A_0>0$ since the solution for $A_0<0$ can be obtained by the transformation: $(t,A)\mapsto (-t,-A)$.
\paragraph{Time asymptotic for \hyperref[par:case1a]{Case \textrm{I}}.} The initial data in this case satisfies $A_0^{\prime} \in (-\infty,-4)$ and $c_0\in\R$. According to \hyperref[item:reg1-i]{\textbf{Regime 1(i)}} and \hyperref[item:reg2-i]{\textbf{Regime 2(i)}}, if we set 
\begin{equation*}
	T_{+}\vcentcolon= \mathcal{I}_{H}(0) =A_0 \int_{1}^{0} \dfrac{\dif s}{W(s)-2}>0, 
\end{equation*}
then $(A(t),A^{\prime}(t)) \to (0,-\infty)$ as $t\nearrow T_{+}$. It follows from (\ref{AIntH}) that
\begin{equation}\label{t-T+}
A_0 \int_0^{A/A_0} \dfrac{\dif s}{W(s)-2} = t - T_+.
\end{equation}
Since $A_0^{\prime}<-4$, (\ref{Daf}) implies that $\beta(A_0^{\prime})<0$. By the expression (\ref{W}), the following series expansion holds:
\begin{equation*}
	\dfrac{1}{W(s)-2} = \Big(\dfrac{|\beta(A_0^{\prime})|}{16}\Big)^{\frac{1}{3}} s^2 + \mathcal{O}(s^4) \quad \text{as } \ s\to 0^{+}.
\end{equation*}
Substituting the above into (\ref{t-T+}), it follows that as $A\to 0^+$,
\begin{equation*}
 t - T_{+}=A_0\int_{0}^{A/A_0}\!\! \Big\{ \Big(\dfrac{|\beta(A_0^{\prime})|}{16}\Big)^{\frac{1}{3}} s^2 + \mathcal{O}(s^4) \Big\}\, \dif s = \dfrac{1}{3A_0^2}\Big(\dfrac{|\beta(A_0^{\prime})|}{16}\Big)^{\frac{1}{3}} A^3 + \mathcal{O}(A^5).
\end{equation*}
Therefore we obtain the convergence rate 
\begin{equation*}
	A(t) = \mathcal{O}\big(\sqrt[\leftroot{-1}\uproot{2}\scriptstyle 3]{T_+-t}\big) \qquad \text{as } \ t\nearrow T_{+}.
\end{equation*}
Combining the above convergence rate with (\ref{W}), we conclude that
\begin{equation*}
A^{\prime}(t) = W\big(\dfrac{A(t)}{A_0}\big) -2 = \mathcal{O}\Big( (T_{+}-t)^{-2/3} \Big) \qquad \text{as } \ t\nearrow T_{+}.
\end{equation*}
Next, we consider the asymptotic of $(A,A^{\prime})$ in the limit $t\to-\infty$. According to the analysis given in \hyperref[item:reg1-ii]{\textbf{Regime 1(ii)}} and \hyperref[item:reg2-ii]{\textbf{Regime 2(ii)}}, as $t\to-\infty$, $A(t)$ takes the form:
\begin{gather}
	A(t)= A_0 \mathcal{I}_2^{-1}\big(\dfrac{t}{A_0}\big), \quad \text{ with } \ \mathcal{I}_2(x)\vcentcolon= \int_{1}^{x} \dfrac{\dif s}{V_2(s)-2}, \nonumber\\
	\text{where } \ V_2(s) = 4 \cos\Big\{ \dfrac{1}{3}\arccos\Big( 1- \dfrac{2}{s^6 \beta(A_0^{\prime})}  \Big) -\dfrac{4\pi }{3} \Big\}. \label{v2}
\end{gather}
It can be verified that $\lim_{s\to\infty}V_2(s)=-2$. This implies $\lim_{x\to\infty}\mathcal{I}_2(x) = -\infty$. Hence its inverse satisfies $\lim_{t\to-\infty}\mathcal{I}_2^{-1}(t) = \infty$. In addition, there exists $N\in\mathbb{N}$ such that for all $x\ge N$, $-5 \le V_2(x)-2\le -3$, hence
\begin{equation*}
	|\mathcal{I}_2(x)| \ge \int_{N}^{x} \dfrac{\dif s}{|V_2(s)-2|} - |\mathcal{I}_2(N) | \ge  \dfrac{x}{5} - \dfrac{N}{5} - |\mathcal{I}_2(N)|, \quad \text{for all } \ x\ge N. 
\end{equation*}
Thus we have $\liminf_{x\to\infty} \frac{|\mathcal{I}_2(x)|}{x} \ge \frac{1}{5}$. For inverse function $\mathcal{I}_2^{-1}(t)$, we have that
\begin{equation*}
	\limsup\limits_{t\to-\infty}\dfrac{|\mathcal{I}_2^{-1}(t)|}{|t|} = \limsup\limits_{x\to\infty}\dfrac{x}{|\mathcal{I}_2(x)|} = \Big(\liminf\limits_{x\to\infty}\dfrac{|\mathcal{I}_2(x)|}{x}\Big)^{-1} \le 5.
\end{equation*}
Therefore we conclude
\begin{equation*}
	A(t)=A_0\mathcal{I}_2^{-1}\big(\dfrac{t}{A_0}\big)=\mathcal{O}(t) \qquad \text{as } \ t\to -\infty.
\end{equation*}
Since $\lim_{x\to\infty}\mathcal{I}_2^{\prime}(x)=-\frac{1}{4}$, by inverse function theorem,
\begin{equation*}
	\lim\limits_{t\to-\infty}A^{\prime}(t) = \lim\limits_{t\to-\infty}\dfrac{1}{\mathcal{I}_2^{\prime}\big( \mathcal{I}_2^{-1}(t/A_0) \big)} = -4, \ \text{ which implies } \ A^{\prime}(t) + 4 = o(1).
\end{equation*}
The following series expansion holds:
\begin{align*}
	&\arccos(1-x^{-6}) = \sqrt{2}x^{-3} + \dfrac{1}{6\sqrt{2}}x^{-9} + \mathcal{O}(x^{-15}) && \text{as } \ x\to \infty,\\
	&4\cos\big( \dfrac{x}{3} - \dfrac{4\pi}{3} \big) + 2 = -\dfrac{2}{\sqrt{3}}x + \dfrac{1}{9}x^2 + \mathcal{O}(x^3) && \text{as } \ x\to 0.
\end{align*}
Using these in the expression (\ref{v2}), we obtain the series expansion:
\begin{equation*}
	V_2(s)+2 = -\Big(\dfrac{16}{3\beta(A_0^{\prime})}\Big)^{\frac{1}{2}} s^{-3} + \dfrac{4}{9\beta(A_0^{\prime})}s^{-6} + \mathcal{O}(s^{-9}) = \mathcal{O}(s^{-3}) \qquad \text{as } \ s\to\infty.
\end{equation*}
Applying the fact that $\mathcal{I}_2^{-1}(t)=\mathcal{O}(t)$ as $t\to-\infty$, we conclude that
\begin{align*}
	&A^{\prime}(t)+4 
	= V_2\big(\mathcal{I}_2^{-1}(\tfrac{t}{A_0})\big) +2 = \mathcal{O}\big(|\mathcal{I}_2^{-1}(\tfrac{t}{A_0})|^{-3}\big) = \mathcal{O}(t^{-3}) \qquad \text{as } \ t\to-\infty.
\end{align*}
In summary if $A_0^{\prime}\in (-\infty,-4)$ and $c_0\in\R$, then the corresponding solution $A(t)$ satisfies the following asymptotic behaviour: there exists $T_{+}=T_{+}(A_0,A_0^{\prime})\in(0,\infty)$ such that
\begin{alignat*}{3}
	&A(t) = \mathcal{O}\big(\sqrt[\leftroot{-1}\uproot{2}\scriptstyle 3]{T_+-t}\big), &&\quad \text{and} \quad A^{\prime}(t) = \mathcal{O}\Big( (T_{+}-t)^{-2/3} \Big) \qquad  && \text{as } \ t\nearrow T_{+},\\
	&A(t)=\mathcal{O}(t), &&\quad \text{and} \quad A^{\prime}(t) = -4 + \mathcal{O}(t^{-3}) && \text{as } \ t\to-\infty.
\end{alignat*}
\paragraph{Time asymptotic for \hyperref[par:case2]{Case \textrm{II}}.} The initial data in this case satisfies $A_0^{\prime} \in (-4,2)$ and $c_0\ge 0$. According to the analysis given in \hyperref[item:reg3-i]{\textbf{Regime 3(i)}} and \hyperref[item:reg4-i]{\textbf{Regime 4(i)}}, in the limit $t\to \infty$, the solution $A(t)$ is takes the form:
\begin{gather}
	A(t)= A_0 \mathcal{I}_0^{-1}\big(\dfrac{t}{A_0}\big), \quad \text{ with } \ \mathcal{I}_0(x)\vcentcolon= \int_{1}^{x} \dfrac{\dif s}{V_0(s)-2}, \nonumber\\
	\text{where } \ V_0(s) = 4 \cos\Big\{ \dfrac{1}{3}\arccos\Big( 1- \dfrac{2}{s^6 \beta(A_0^{\prime})}  \Big) \Big\}. \label{v0}
\end{gather}
It can be verified that $\lim_{s\to\infty}V_0(s)=4$. This implies $\lim_{x\to\infty}\mathcal{I}_0(x) = \infty$. Hence its inverse satisfies $\lim_{t\to\infty}\mathcal{I}_0^{-1}(t) = \infty$. In addition, there exists $N\in\mathbb{N}$ such that for all $x\ge N$, $1 \le V_0(x)-2\le 3$, hence
\begin{equation*}
	|\mathcal{I}_0(x)| \ge \int_{N}^{x} \dfrac{\dif s}{|V_0(s)-2|} - |\mathcal{I}_0(N) | \ge  \dfrac{x}{3} - \dfrac{N}{3} - |\mathcal{I}_0(N)|, \quad \text{for all } \ x\ge N. 
\end{equation*}
Thus we have $\liminf_{x\to\infty} \frac{|\mathcal{I}_0(x)|}{x} \ge \frac{1}{3}$. For inverse function $\mathcal{I}_0^{-1}(t)$, we have that
\begin{equation*}
	\limsup\limits_{t\to\infty}\dfrac{|\mathcal{I}_0^{-1}(t)|}{t} = \limsup\limits_{x\to\infty}\dfrac{x}{|\mathcal{I}_0(x)|} = \Big(\liminf\limits_{x\to\infty}\dfrac{|\mathcal{I}_0(x)|}{x}\Big)^{-1} \le 3.
\end{equation*}
Therefore we conclude
\begin{equation*}
	A(t)=A_0\mathcal{I}_0^{-1}\big(\dfrac{t}{A_0}\big)=\mathcal{O}(t) \qquad \text{as } \ t\to \infty.
\end{equation*}
Since $\lim_{x\to\infty}\mathcal{I}_0^{\prime}(x)=\frac{1}{2}$, by inverse function theorem,
\begin{equation*}
	\lim\limits_{t\to\infty}A^{\prime}(t) = \lim\limits_{t\to\infty}\dfrac{1}{\mathcal{I}_0^{\prime}\big( \mathcal{I}_0^{-1}(t/A_0) \big)} = 2, \ \text{ which implies } \ A^{\prime}(t) - 2 = o(1).
\end{equation*}
The following series expansion holds:
\begin{align*}
	&\arccos(1-x^{-6}) = \sqrt{2}x^{-3} + \dfrac{1}{6\sqrt{2}}x^{-9} + \mathcal{O}(x^{-15}) && \text{as } \ x\to \infty,\\
	&4\cos\big( \dfrac{x}{3} \big) - 4 = -\dfrac{2}{9}x^2 + \dfrac{1}{486}x^4 + \mathcal{O}(x^5) && \text{as } \ x\to 0.
\end{align*}
Using these in the expression (\ref{v0}), we obtain the series expansion:
\begin{equation*}
	V_0(s)-4 = -\dfrac{8}{9\beta(A_0^{\prime})} s^{-6} + \dfrac{8}{243\beta(A_0^{\prime})}s^{-12} + \mathcal{O}(s^{-15}) = \mathcal{O}(s^{-6}) \qquad \text{as } \ s\to\infty.
\end{equation*}
Applying the fact that $\mathcal{I}_0^{-1}(t)=\mathcal{O}(t)$ as $t\to \infty$, we conclude that
\begin{align*}
	&A^{\prime}(t)-2 = V_0\big(\mathcal{I}_0^{-1}(\tfrac{t}{A_0})\big) -4 = \mathcal{O}\big(|\mathcal{I}_0^{-1}(\tfrac{t}{A_0})|^{-6}\big) = \mathcal{O}(t^{-6}) \qquad \text{as } \ t\to \infty.
\end{align*}
Next we consider the asymptotic behaviour of $A(t)$, $A^{\prime}(t)$ in the limit $t\to-\infty$. In this case, according to \hyperref[item:reg3-ii]{\textbf{Regime 3(ii)}} and \hyperref[item:reg4-ii]{\textbf{Regime 4(ii)}}, the solution $A(t)$ is given by
\begin{gather}
	A(t)= A_0 \mathcal{I}_1^{-1}\big(\dfrac{t}{A_0}\big), \quad \text{ with } \ \mathcal{I}_1(x)\vcentcolon= \int_{1}^{x} \dfrac{\dif s}{V_1(s)-2}, \nonumber\\
	\text{where } \ V_1(s) = 4 \cos\Big\{ \dfrac{1}{3}\arccos\Big( 1- \dfrac{2}{s^6 \beta(A_0^{\prime})}  \Big) - \dfrac{2\pi}{3} \Big\}. \label{v1}
\end{gather}
It can be verified that $\lim_{s\to\infty}V_1(s)=-2$. This implies $\lim_{x\to\infty}\mathcal{I}_1(x) = -\infty$. Hence its inverse satisfies $\lim_{t\to-\infty}\mathcal{I}_1^{-1}(t) = \infty$. In addition, there exists $N\in\mathbb{N}$ such that for all $x\ge N$, $-5 \le V_0(x)-2\le -3$, hence
\begin{equation*}
	|\mathcal{I}_1(x)| \ge \int_{N}^{x} \dfrac{\dif s}{|V_1(s)-2|} - |\mathcal{I}_1(N) | \ge  \dfrac{x}{5} - \dfrac{N}{5} - |\mathcal{I}_1(N)|, \quad \text{for all } \ x\ge N. 
\end{equation*}
Thus we have $\liminf_{x\to\infty} \frac{|\mathcal{I}_1(x)|}{x} \ge \frac{1}{5}$. For inverse function $\mathcal{I}_1^{-1}(t)$, we have that
\begin{equation*}
	\limsup\limits_{t\to-\infty}\dfrac{|\mathcal{I}_1^{-1}(t)|}{|t|} = \limsup\limits_{x\to\infty}\dfrac{x}{|\mathcal{I}_1(x)|} = \Big(\liminf\limits_{x\to\infty}\dfrac{|\mathcal{I}_1(x)|}{x}\Big)^{-1} \le 5.
\end{equation*}
Therefore we conclude
\begin{equation*}
	A(t)=A_0\mathcal{I}_1^{-1}\big(\dfrac{t}{A_0}\big)=\mathcal{O}(t) \qquad \text{as } \ t\to -\infty.
\end{equation*}
Since $\lim_{x\to\infty}\mathcal{I}_1^{\prime}(x)=-\frac{1}{4}$, by inverse function theorem,
\begin{equation*}
	\lim\limits_{t\to-\infty}A^{\prime}(t) = \lim\limits_{t\to-\infty}\dfrac{1}{\mathcal{I}_1^{\prime}\big( \mathcal{I}_1^{-1}(t/A_0) \big)} = -4, \ \text{ which implies } \ A^{\prime}(t) + 4 = o(1).
\end{equation*}
The following series expansion holds:
\begin{align*}
	&\arccos(1-x^{-6}) = \sqrt{2}x^{-3} + \dfrac{1}{6\sqrt{2}}x^{-9} + \mathcal{O}(x^{-15}) && \text{as } \ x\to \infty,\\
	&4\cos\big( \dfrac{x}{3} - \dfrac{2\pi}{3} \big) + 2 = \dfrac{2}{\sqrt{3}}x + \dfrac{1}{9}x^2 + \mathcal{O}(x^3) && \text{as } \ x\to 0.
\end{align*}
Using these in the expression (\ref{v1}), we obtain the series expansion:
\begin{equation*}
	V_1(s)+2 = \Big(\dfrac{16}{3\beta(A_0^{\prime})}\Big)^{\frac{1}{2}} s^{-3} + \dfrac{4}{9\beta(A_0^{\prime})}s^{-6} + \mathcal{O}(s^{-9}) = \mathcal{O}(s^{-3}) \qquad \text{as } \ s\to\infty.
\end{equation*}
Applying the fact that $\mathcal{I}_1^{-1}(t)=\mathcal{O}(t)$ as $t\to-\infty$, we conclude that
\begin{align*}
	&A^{\prime}(t)+4 
	= V_1\big(\mathcal{I}_1^{-1}(\tfrac{t}{A_0})\big) +2 = \mathcal{O}\big(|\mathcal{I}_1^{-1}(\tfrac{t}{A_0})|^{-3}\big) = \mathcal{O}(t^{-3}) \qquad \text{as } \ t\to-\infty.
\end{align*}
In summary, if $A_0^{\prime}\in(-4,2)$ and $c_0\ge 0$, then the corresponding solution $A(t)$ satisfies the following time asymptotic behaviour:
\begin{align*}
&A(t) = \mathcal{O}(t), \quad \text{and} \quad  A^{\prime}(t)=2+\mathcal{O}(t^{-6}) && \text{as } \ t\to \infty,\\
&A(t) = \mathcal{O}(t), \quad \text{and} \quad A^{\prime}(t)=-4+\mathcal{O}(t^{-3}) && \text{as } \ t\to -\infty.
\end{align*}
\paragraph{Time asymptotic for \hyperref[par:case3a]{ Case \textrm{III}}.} The initial data in this case satisfies $A_0^{\prime} \in (2,\infty)$ and $c_0\in\R$. According to \hyperref[item:reg5-i]{\textbf{Regime 5(i)}}, as $t\to\infty$, the solution $A(t)$ takes the form:
\begin{equation*}
	A(t)= A_0 \mathcal{I}_H^{-1}\big(\dfrac{t}{A_0}\big), \quad \text{ with } \ \mathcal{I}_H(x)\vcentcolon= \int_{1}^{x} \dfrac{\dif s}{W(s)-2},
\end{equation*}
where $W(s)$ is defined in (\ref{W}). It can be verified that $\lim_{s\to\infty}W(s)=4$. This implies $\lim_{x\to\infty}\mathcal{I}_H(x) = \infty$. Hence its inverse satisfies $\lim_{t\to\infty}\mathcal{I}_H^{-1}(t) = \infty$. In addition, there exists $N\in\mathbb{N}$ such that for all $x\ge N$, $1 \le W(x)-2\le 3$, hence
\begin{equation*}
	|\mathcal{I}_H(x)| \ge \int_{N}^{x} \dfrac{\dif s}{|W(s)-2|} - |\mathcal{I}_H(N) | \ge  \dfrac{x}{3} - \dfrac{N}{3} - |\mathcal{I}_H(N)|, \quad \text{for all } \ x\ge N. 
\end{equation*}
Thus we have $\liminf_{x\to\infty} \frac{|\mathcal{I}_H(x)|}{x} \ge \frac{1}{3}$. For inverse function $\mathcal{I}_H^{-1}(t)$, we have that
\begin{equation*}
	\limsup\limits_{t\to\infty}\dfrac{|\mathcal{I}_H^{-1}(t)|}{t} = \limsup\limits_{x\to\infty}\dfrac{x}{|\mathcal{I}_H(x)|} = \Big(\liminf\limits_{x\to\infty}\dfrac{|\mathcal{I}_H(x)|}{x}\Big)^{-1} \le 3.
\end{equation*}
Therefore we conclude
\begin{equation*}
	A(t)=A_0\mathcal{I}_H^{-1}\big(\dfrac{t}{A_0}\big)=\mathcal{O}(t) \qquad \text{as } \ t\to \infty.
\end{equation*}
Since $\lim_{x\to\infty}\mathcal{I}_H^{\prime}(x)=\frac{1}{2}$, by inverse function theorem,
\begin{equation*}
	\lim\limits_{t\to \infty}A^{\prime}(t) = \lim\limits_{t\to \infty}\dfrac{1}{\mathcal{I}_H^{\prime}\big( \mathcal{I}_H^{-1}(t/A_0) \big)} = 2, \ \text{ which implies } \ A^{\prime}(t) -2  = o(1).
\end{equation*}
From (\ref{W}), one has the following series expansion:
\begin{equation*}
W(s) -4  = -\dfrac{8(3|\beta(A_0^{\prime})|^2-2)}{9(\beta(A_0^{\prime}))^3} s^{-6} + \mathcal{O}(s^{-9}) = \mathcal{O}(s^{-6}) \qquad \text{as } \ s\to\infty.
\end{equation*}
Applying the fact that $\mathcal{I}_H^{-1}(t)=\mathcal{O}(t)$ as $t\to\infty$, we conclude that
\begin{align*}
	&A^{\prime}(t)-2  
	= W\big(\mathcal{I}_H^{-1}(\tfrac{t}{A_0})\big) -4 = \mathcal{O}\big(|\mathcal{I}_H^{-1}(\tfrac{t}{A_0})|^{-6}\big) = \mathcal{O}(t^{-6}) \qquad \text{as } \ t\to\infty.
\end{align*}
Next, we consider the asymptotic of $(A,A^{\prime})$ as $t$ decreases. According to the analysis given in \hyperref[item:reg5-ii]{\textbf{Regime 5(ii)}}, if we set 
\begin{equation*}
	T_{-}\vcentcolon= \mathcal{I}_{H}(0) =A_0 \int_{1}^{0} \dfrac{\dif s}{W(s)-2}<0, 
\end{equation*}
then $(A(t),A^{\prime}(t)) \to (0,\infty)$ as $t\searrow T_{-}$. It follows from (\ref{AIntH}) that
\begin{equation}\label{t-T-}
	A_0 \int_0^{A/A_0} \dfrac{\dif s}{W(s)-2} = t - T_{-}.
\end{equation}
Since $A_0^{\prime}>2$, (\ref{Daf}) implies that $\beta(A_0^{\prime})>0$. By the expression (\ref{W}), the following series expansion holds:
\begin{equation*}
	\dfrac{1}{W(s)-2} = -\Big(\dfrac{|\beta(A_0^{\prime})|}{16}\Big)^{\frac{1}{3}} s^2 + \mathcal{O}(s^4) \quad \text{as } \ s\to 0^{+}.
\end{equation*}
Substituting the above into (\ref{t-T-}), it follows that as $A\to 0^+$,
\begin{equation*}
	t - T_{-}=A_0\int_{0}^{A/A_0}\!\! \Big\{- \Big(\dfrac{|\beta(A_0^{\prime})|}{16}\Big)^{\frac{1}{3}} s^2 + \mathcal{O}(s^4) \Big\}\, \dif s = -\dfrac{1}{3A_0^2}\Big(\dfrac{|\beta(A_0^{\prime})|}{16}\Big)^{\frac{1}{3}} A^3 + \mathcal{O}(A^5).
\end{equation*}
Therefore we obtain the convergence rate 
\begin{equation*}
	A(t) = \mathcal{O}\big(\sqrt[\leftroot{-1}\uproot{2}\scriptstyle 3]{t-T_{-}}\big) \qquad \text{as } \ t\searrow T_{-}.
\end{equation*}
Combining the above convergence rate with (\ref{W}), we conclude that
\begin{equation*}
	A^{\prime}(t) = W\big(\dfrac{A(t)}{A_0}\big) -2 = \mathcal{O}\Big( (t-T_{-})^{-2/3} \Big) \qquad \text{as } \ t\searrow T_{-}.
\end{equation*}
In summary if $A_0^{\prime}\in (2,\infty)$ and $c_0\in\R$, then the corresponding solution $A(t)$ satisfies the following asymptotic behaviour: there exists $T_{-}=T_{-}(A_0,A_0^{\prime})\in(-\infty,0)$ such that
\begin{alignat*}{3}
	&A(t)=\mathcal{O}(t), &&\quad \text{and} \quad A^{\prime}(t) = 2 + \mathcal{O}(t^{-6}) \qquad &&\text{as } \ t\to\infty,\\
	&A(t)=\mathcal{O}\big(\sqrt[\leftroot{-1}\uproot{2}\scriptstyle 3]{t-T_{-}}\big), &&\quad \text{and} \quad A^{\prime}(t) = \mathcal{O}\Big( (t-T_{-})^{-2/3} \Big) \qquad &&\text{as } \ t\searrow T_{-}.
\end{alignat*}

\section{Flow of Parabola with Drift}\label{sec:4}
Suppose that $A(t)$ is a solution to (\ref{A''})--(\ref{A}). We consider the following ansatz:
\begin{equation}\label{drift}
	\phi(t,x,y,z) = \dfrac{x^2+y^2-2z^2}{A(t)} + B(t) z 
\end{equation}
For each point $\x\in\d\Omega(t)$ at time $t$, let $\X(s;t,\x)=(X^x,X^y,X^z)(s;t,\x)$ be the backward characteristic curve emanating from $(t,\x)$, which satisfies $\X(t;t,x)=\x$ and the differential equations
\begin{equation*}
	\dfrac{\d X^x}{\d s}=\dfrac{2 X^x}{A(s)}, \qquad \dfrac{\d X^y}{\d s}=\dfrac{2 X^y}{A(s)}, \qquad 
	\dfrac{\d X^z}{\d s} = B(s)- \dfrac{4 X^z}{A(s)}.
\end{equation*}
For simplicity denote $X^\sigma_0\equiv X^\sigma(0;t,\x) $ for $\sigma=x,\,y,\,z$. Solving these, one has
\begin{gather}
X^x_0=x e^{-g(t)}, \quad X^y_0 = y e^{-g(t)},\ \quad X^z_0= z e^{2g(t)} - h(t),\label{dc}\\
\text{where } \quad h(t)\vcentcolon= \int_{0}^{t}\!\! B(s) e^{2g(s)}\, \dif s, \quad \text{ and } \quad g(t)=\int_{0}^{t}\!\! \dfrac{2}{A(s)}\, \dif s.\nonumber
\end{gather}
Therefore, the equation describing the free surface must satisfies the form:
\begin{equation}\label{gSd}
	f\big(X^x_0,X^y_0,X^z_0\big) =f\big(x e^{-g(t)},y e^{-g(t)},z e^{2g(t)} - h(t)\big)=0,
\end{equation}
for some function $f(\cdot,\cdot,\cdot)$ independent of $(t,x,y,z)$. On the other hand, substituting (\ref{drift}) into (\ref{bernoulli}), we have
\begin{equation}\label{gBd}
	\Big\{ \dfrac{2-A^{\prime}}{A^2} (x^2+y^2) + \dfrac{2(4+A^{\prime})}{A^2} z^2 + \big( B^{\prime} - \dfrac{4B}{A} \big) z \Big\} \Big\vert_{F(t,\x)=0} = c(t) - \dfrac{B^2}{2}.
\end{equation}
Rewriting the above in terms of $\big(X^x_0,X^y_0,X^z_0\big)$ given in (\ref{dc}), one obtains
\begin{align*}
	0=& \dfrac{2-A^{\prime}}{A^2} e^{2g} |X^x_0|^2+ \dfrac{2-A^{\prime}}{A^2} e^{2g} |X^y_0|^2 +\dfrac{2(4+A^{\prime})}{A^2}e^{-4g} |X^z_0|^2\\
	&+ e^{-2g} \Big( \dfrac{4(4+A^{\prime})}{A^2} h e^{-2g} + B^{\prime} - \dfrac{4B}{A} \Big) X^z_0 \\
	&+ h e^{-2g} \Big( B^{\prime} - \dfrac{4B}{A} + \dfrac{2(4+A^{\prime})}{A^2}he^{-2g} \Big) + \dfrac{B^2}{2} - c(t).
\end{align*}
Since the above equation must be consistent with (\ref{gSd}), it follows that there exists constants $c_1$, $c_2$, $c_3\in\R$ such that
\begin{subequations}\label{consisd}
	\begin{gather}
		\dfrac{2-A^{\prime}}{A^2} e^{2g} = c_1 \dfrac{2(4+A^{\prime})}{A^2} e^{-4g},\label{consisd1}\\
		\dfrac{2-A^{\prime}}{A^2} e^{2g} = c_2 e^{-2g} \Big( \dfrac{4(4+A^{\prime})}{A^2} h e^{-2g} + B^{\prime} - \dfrac{4B}{A} \Big),\label{consisd2}\\
		\dfrac{2-A^{\prime}}{A^2} e^{2g} =c_3 \Big\{ \big( B^{\prime} - \dfrac{4B}{A}\big) h e^{-2g} + \dfrac{2(4+A^{\prime})}{A^2}h^2e^{-4g}  + \dfrac{B^2}{2} - c(t) \Big\}.\label{consisd3}
	\end{gather}
\end{subequations}
From (\ref{consisd1}), we have the identity $\frac{1}{2c_1}=\frac{4+A^{\prime}}{2-A^{\prime}}e^{-6g}$. Substituting this into (\ref{consisd2}),
\begin{align*}
	1=\dfrac{2c_2}{c_1} h + c_2 \Big( B^{\prime} A - 4B \Big) \dfrac{A}{2-A^{\prime}}e^{-4g} \ \iff \ h=\dfrac{c_1}{2} \Big\{ \dfrac{1}{c_2} + \dfrac{4 AB - B^{\prime}A^2}{2-A^{\prime}} e^{-4g} \Big\}.
\end{align*}
Taking derivative on the above equation, and using the definition $h(t)=\int_{0}^{t}\!\!B(s)e^{2g(s)}\,\dif s$, we obtain 
\begin{equation*}
	0=\dfrac{2B}{c_1}+ \Big\{ B^{\prime\prime}A + 4(A^{\prime}-1)B^{\prime} - \dfrac{12A^{\prime}}{A}B  \Big\}\dfrac{A}{2-A^{\prime}}e^{-6g}.
\end{equation*}
Using the identity $\frac{1}{2c_1}=\frac{4+A^{\prime}}{2-A^{\prime}}e^{-6g}$ on the above equation once again, then rearranging terms, we have
\begin{equation}\label{B}
	B^{\prime\prime} + \dfrac{4(A^{\prime}-1)}{A} B^{\prime} + \dfrac{8(2-A^{\prime})}{A^2} B = 0.
\end{equation}
Assuming further that $A(t)=-4t+4a_0$ where $a_0=\frac{A_0}{4}\neq 0$, then one has the equation
\begin{gather*}
	B^{\prime\prime} + \dfrac{5B^{\prime}}{t-a_0}  + \dfrac{3 B}{(t-a_0)^2} = 0, \quad \text{ which implies } \quad B(t) = \dfrac{k_1}{t-a_0} + \dfrac{k_2}{(t-a_0)^3},\\
	\text{where } \quad k_1\vcentcolon=\frac{a_0}{2}(a_0 B_0^{\prime}-3B_0), \quad k_2\vcentcolon=\frac{a_0^3}{2}(B_0-a_0 B_0^{\prime}), \quad (B_0,B_0^{\prime})\vcentcolon=(B(0),B^{\prime}(0)).
\end{gather*}
In addition, using the expression (\ref{consisd3}), one can calculate to get that $$c(t)-\frac{B^2}{2} = \Big( c_0 - \dfrac{B_0^2}{2} \Big)\frac{a_0^3}{(a_0-t)^3}.$$ where $c_0\equiv c(0)$. Substituting these expression of $A$ and $B$ into (\ref{gBd}), we conclude that the free surface is a \textbf{Paraboloid} described by the equation
\begin{equation}\label{paraboloid}
	z=\dfrac{3(t-a_0)^2}{16 k_2} ( x^2 + y^2 ) + \dfrac{a_0^3}{2k_2} \Big(c_0-\dfrac{B_0^2}{2}\Big)(t-a_0).
\end{equation}
As $|t|\to\infty$, this paraboloid will collapse to \textbf{$z$-Axis} described by $x^2+y^2=0$.


\section{Flow of Ellipsoid with Constant Vorticity}\label{sec:5}
In this section we construct an explicit solution to (\ref{lagrange})--(\ref{bernoulli}) with vorticity, i.e. when the condition (\ref{harmonic}) is not satisfied. More precisely, the solution has constant vorticity $ \curl \vv = -2\ep \hz$ for some constant $\ep>0$. We remark that such solution was first studied by Ovsjannikov in \cite{Ovsyannikov}.

Let $\X(t) = (x(t),y(t),z(t))$ be the position vector of the free surface at time $t\ge 0$, and denote $\X_0 \vcentcolon= \X(0) =\vcentcolon (x_0,y_0,z_0)$. We impose the following ansatz:
\begin{equation}\label{vorAns}
\X(t)=A(t)\cdot \X_0 \qquad \text{where} \qquad A(t)\vcentcolon=
\begin{pmatrix}
	\alpha& \beta& 0\\
	-\beta &\alpha& 0\\
	0&0&\gamma
\end{pmatrix}.
\end{equation}
Here $\alpha(t),\, \beta(t),\, \gamma(t)$ are functions of 
$t$ such that $A(0)=I$ with $I$ being the $3$-by-$3$ identity matrix. Let $P=P(t,\X)$ be the pressure. By chain rule one has,
\begin{equation*}
\nabla_{\X_0} \big(P(t,\X) \big) = \nabla_{\X_0} \big( P(t,A\cdot \X_0) \big) = A^{\top} \cdot \nabla_{\X} P\big\vert_{\X=A\X_0}.
\end{equation*}
Then Newton's second law yields
\begin{equation*}
\X''(t)=-\nabla_{\X}P\big\vert_{\X=A\X_0}=-(A^{\top})^{-1}\cdot \nabla_{\X_0} \big( P(t,\X) \big). 
\end{equation*}
Substituting the ansatz $\X(t)=A(t)\X_0$ into the above equation, we obtain that, 
\begin{equation*}
	A^{\top}(t) \cdot A''(t)  \cdot \X_0= -\nabla_{\X_0} \big( P(t,\X) \big).
\end{equation*}
Integrating in $\X_0$ and noting that on the free surface, pressure is constant for given time $t$, we get that 
\begin{equation*}
	\frac{1}{2} \X_0^{\top} \cdot A^{\top}(t) \cdot A''(t) \cdot \X_0=c(t).
\end{equation*}
In addition, we also impose the following similarity relations:
\begin{equation}\label{SimC}
	A^{\top} A''=a(t)I, \qquad c(t)=a(t) \dfrac{|\X_0|^2}{2},
\end{equation}
Our aim is then to find $\alpha, \beta, \gamma, $ and $a$ that solve the problem with initial conditions:
\begin{equation}
	\alpha(0)=1, \quad \beta(0)=0, \quad \gamma(0)=1, \quad a(0)=a_0.
\end{equation}
Since $\X(t)$ is also the particle trajectory, if we set $\vv(t,\X)$ as the velocity, then
\begin{equation*}
\vv\big(t,\X(t)\big) = \dfrac{\dif \X}{\dif t}(t) = A^{\prime} (t) \X_0 = A^{\prime} A^{-1} \X(t).
\end{equation*}
From the incompressibility $\div\,\vv = 0$, and Jacobi's formula for determinant, one has
\begin{equation*}
	0=\div_{\X} \big( \vv (t,\X) \big)=\textrm{tr}(A^{\prime} A^{-1}) = \dfrac{\dif }{\dif t} \ln(\det A),
\end{equation*}
It follows that $\det A(t) = \det A(0) = 1$, hence by (\ref{vorAns}),
\begin{equation}\label{eq:det-1}
	\alpha^2+\beta^2=\frac1\gamma.
\end{equation}
Next, by a direct computation one has
\begin{equation*}
	A^{\top}A''
	=
	\begin{pmatrix}
		\alpha& -\beta& 0\\
		\beta &\alpha& 0\\
		0&0&\gamma
	\end{pmatrix}
	\begin{pmatrix}
		\alpha''& \beta''& 0\\
		-\beta'' &\alpha''& 0\\
		0&0&\gamma''
	\end{pmatrix}
	=
	\begin{pmatrix}
		\alpha\alpha''+\beta\beta''& \alpha\beta''-\beta\alpha''& 0\\
		\beta\alpha''-\alpha\beta'' &\beta\beta''+\alpha\alpha''& 0\\
		0&0&\gamma\gamma''
	\end{pmatrix}.
\end{equation*}
Then the similarity assumption (\ref{SimC}) yields
\begin{equation}\label{eq:ddnnbbbcc}
\left\{
\begin{array}{lll}
	\alpha\beta''-\alpha''\beta=0,\\
	\alpha\alpha''+\beta\beta''=\gamma\gamma''. 
	\end{array}
	\right.
\end{equation}
From here, we wish to reduce the system into a single ODE for 
$\gamma$ and solve it. Using \eqref{eq:ddnnbbbcc}, we get that 
$(\alpha\beta'-\alpha'\beta)'=0,$ hence
after integration we have 
\begin{equation*}
	\alpha\beta'-\alpha'\beta=\epsilon=\mbox{const}.
\end{equation*}
Moreover, differentiating \eqref{eq:det-1}, we have $\alpha\alpha'+\beta\beta=-\frac{\gamma'}{2\gamma^2}$. Putting these equations together, it forms the system: 
\begin{eqnarray}
	\alpha\beta'-\alpha'\beta=\epsilon, \label{eq:epsilon-yytt} \\
	\alpha\alpha'+\beta\beta=-\frac{\gamma'}{2\gamma^2}.\label{eq:hingtas}
\end{eqnarray}
Multiplying the first equation by $\alpha$, and the second one by $\beta$
and adding up we get 
$\alpha^2\beta'-\alpha\alpha'\beta+\alpha\alpha'\beta+\beta^2\beta'=\alpha\epsilon-\frac{\gamma'\beta}{2\gamma^2}$, 
or 
equivalently
\begin{equation*}
	(\alpha^2+\beta^2)\beta'=\alpha\epsilon-\frac{\gamma'\beta}{2\gamma^2}.
\end{equation*}
Utilizing \eqref{eq:det-1} we obtain
\begin{equation*}
	\beta'=\alpha\epsilon\gamma-\frac{\gamma'\beta}{2\gamma}.
\end{equation*}
Repeating the same calculation but with $\alpha$ and $\beta$ exchanged, we also obtain 
\begin{equation*}
	\alpha'=-\beta\epsilon\gamma-\frac{\gamma'\alpha}{2\gamma}.
\end{equation*}
Using the last two equations and \eqref{eq:det-1}
one more time we get that 
\begin{equation}\label{a'b'}
	(\alpha')^2+(\beta')^2=\alpha^2\epsilon^2\gamma^2+\frac{(\gamma')^2\beta^2}{4\gamma^2}+
	\beta^2\epsilon^2\gamma^2+
	\frac{(\gamma')^2\alpha^2}{4\gamma^2}
	=
	\epsilon^2\gamma
	+
	\frac{(\gamma')^2}{4\gamma^3}.
\end{equation}
Differentiating \eqref{eq:hingtas} yields
\[
(\alpha')^2+\alpha\alpha''+(\beta')^2
+\beta\beta''=-\frac{\gamma''}{2\gamma^2}+
\frac{(\gamma')^2}{\gamma^3}.
\]
Substituting (\ref{eq:ddnnbbbcc}) and (\ref{a'b'}) into the above, we obtain
\[
\epsilon^2\gamma+\frac{(\gamma')^2}{4\gamma^3}+
\gamma\gamma''=-\frac{\gamma''}{2\gamma^2}+
\frac{(\gamma')^2}{\gamma^3}.
\]
Hence we get the desired ODE for 
$\gamma$
\begin{equation}\label{eq:pepperton-00}
	(\gamma^3+\frac12)\gamma\gamma''-\frac34(\gamma')^2+\epsilon^2\gamma^4=0, \qquad \gamma(0)=1.
\end{equation}

\begin{figure}
\includegraphics[trim={1.2cm 0 0 0},clip,scale=0.7]{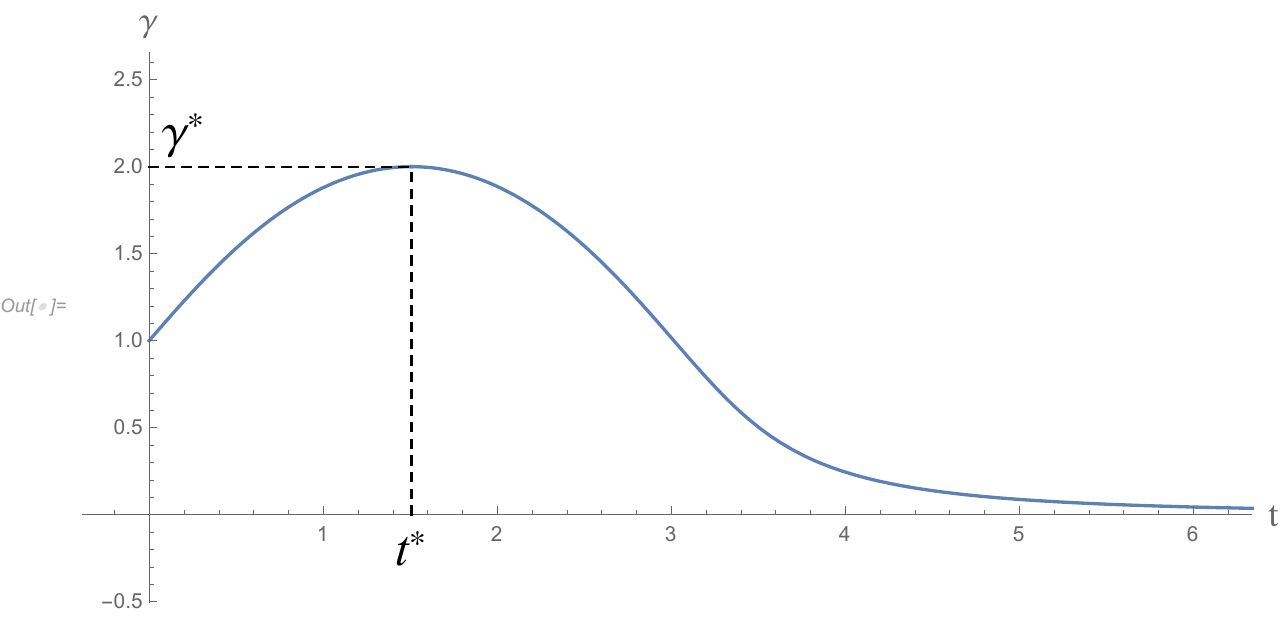}
\caption{The plot of $\gamma$.}
\label{fig-ovs}
\end{figure}

In order to solve (\ref{eq:pepperton-00}) we make the 
substitution $2\psi(\gamma)=(\gamma')^2$, 
for a new unknown $\psi$. Then 
$\psi $ solves the following first order ODE
\[
(\gamma^3+\frac12)\gamma\psi'(\gamma)-\frac34\psi(\gamma)+\epsilon^2\gamma^4=0.
\]
Note that 
\[
\frac{d}{d\gamma}
\left(
\frac{(\gamma^3+\frac12)\psi'(\gamma)}{\gamma^3}
\right)
=
\frac{d}{d\gamma}
\left((1+\frac1{2\gamma^3})\psi(\gamma)
\right)
=-\frac3{2\gamma^4}
+
\frac{\gamma^3+\frac12}{\gamma^3}\psi'(\gamma).
\]
This yields, 
\[
\frac{d}{d\gamma}
\left(
\frac{(\gamma^3+\frac12)\psi(\gamma)}{\gamma^3}
+\epsilon^2\gamma
\right)=0
\quad \text{ or } \quad
\frac{(\gamma^3+\frac12)\psi(\gamma)}{\gamma^3}
+\epsilon^2\gamma=C.
\]
Recalling the definition of 
$\psi$, we finally get 
\begin{equation*}
	(\gamma')^2=\frac{2\gamma^3(C-\epsilon^2\gamma)}{\gamma^3+\frac12}.
\end{equation*}
Denote 
$\gamma^*=C/\epsilon^2$, then 
\begin{equation}\label{eq:vetzovs}
	(\gamma')^2=\frac{2\epsilon^2\gamma^3(\gamma^*-\gamma)}{\gamma^3+\frac12}.
\end{equation}
We immediately see that 
$\gamma\le \gamma^*$. 
Using the initial condition 
$\gamma(0)=1$, we get that 
\[
\gamma'(0)=2|\epsilon|\sqrt{\frac{\gamma^*-1}{3}}, \quad \gamma^*\ge 1.
\]
For this case, if $t>0$ and close to $t=0$, 
$\gamma(t)$ must be increasing. Thus integrating 
\eqref{eq:vetzovs}, we have
\[
\int_1^\gamma
\gamma^{-\frac32}\sqrt{2\gamma^3+1}\frac{d\gamma}{\sqrt{\gamma^*-\gamma}}
=
2|\epsilon| t.
\]
As $\gamma$ increases, it attains the value $\gamma^*$
at $t=t^*$, where $t^*$
is determined from 
\[
\int_1^{\gamma^*}
\gamma^{-\frac32}\sqrt{2\gamma^3+1}\frac{d\gamma}{\sqrt{\gamma^*-\gamma}}
=
2|\epsilon| t^*.
\]
From \eqref{eq:pepperton-00}
and \eqref{eq:vetzovs}
it follows that $\gamma$ has a local 
maximum at $t=t^{*}$, and 
$\gamma(t)$
decrease for $t>t^*$. For these values of $t$ we have 
\[
\int^{\gamma^*}_\gamma
\gamma^{-\frac32}\sqrt{2\gamma^3+1}\frac{d\gamma}{\sqrt{\gamma^*-\gamma}}
=
2|\epsilon| (t-t^*).
\]
From here it follows that as $t\to \infty$
then $\gamma\to 0$. The graph of 
$\gamma$ is shown in Figure \ref{fig-ovs}.

The functions $\alpha, \beta$ can be calculated as functions of $\gamma$ as follows: in light of the incompressibility condition \eqref{eq:det-1}, and constant vorticity condition (\ref{eq:epsilon-yytt}), we let 
\[
\alpha=\frac1{\sqrt\gamma}\cos\theta, 
\quad 
\beta=\frac1{\sqrt\gamma}\sin\theta.
\]
To determine $\theta$,
we use \eqref{eq:epsilon-yytt} to obtain that
\[
\alpha
=\frac1{\sqrt\gamma}
\cos 
\left(
\epsilon \int_0^t \gamma dt
\right), 
\quad 
\beta
=\frac1{\sqrt\gamma}
\sin
\left(
\epsilon \int_0^t \gamma dt
\right).
\]

Now let us discuss the 
geometric picture for this flow. The pressure can be assumed 
zero on the sphere 
$\xi^2+\eta^2+\zeta^2=1$, and the 
mapping given by the matrix $A$ implies
\[
\gamma(x^2+y^2)+\frac1{\gamma^2}z^2=1.
\]
We see that 
for $0<t<t^*$ the sphere $\Gamma_0$ at $t=0$ turns into an 
ellipsoid of revolution with axis $z$ until the moment 
when its major semiaxis becomes equal to 
$\gamma^*$. After that the ellipsoid reverses its motion back to the unit sphere, and consequently flattens out to the plane $z=0$, as $t\to \infty$ (see Figure \ref{fig-ovs}). The 3D animation of this flow can be viewed here; 
\url{https://www.maths.ed.ac.uk/~aram/mov.html} .

To clarify the mechanical phenomenon, let us look at the velocity $\vv =A'(0)\X_0$ at initial data. Since $\beta'(0)=\epsilon, \gamma'(0)=2|\epsilon|\sqrt{\frac{\gamma^*-1}3}$, and 
$\alpha'(0)=-\frac12\gamma'(0)$, we have
\[
A^{\prime}(0)=\begin{pmatrix}
	\alpha'(0)&\epsilon&0\\
	-\epsilon&\alpha'(0)&0\\
	0&0&-2\alpha'(0)
\end{pmatrix}.
\]
Observe that the matrix 
$A'(0)$ is not symmetric, which indicates that the motion is not irrotational. In fact $\curl \vv=-2\epsilon \hz$.
Therefore the initial state of the motion, 
at $t=0$, corresponds to uniform rotation of the 
sphere $\Gamma_0$ about $z$ axis 
with angular velocity $\epsilon.$
If $\epsilon=0$ then we get a 
potential flow discussed earlier, for which the needle 
stretches to a line as $t\to\infty$.
However, for $\epsilon\not=0$, the stretching into line 
is impossible and the flow collapses to a hyperplane, as explained above.

The velocity at every point can be written in the form 
$\nabla \phi+\theta$, where 
$\phi$ is a potential function and $\theta$ is the rotation component about the $z$-axis. The kinematic energy of the system is the sum 
of kinetic energies 
\[
K=K_{\nabla\phi}+K_\theta.
\]
During the flow $K_{\nabla\phi}$ transfers to 
$K_{\theta}$ up to the point $t=t^*$, where
$K_{\nabla\phi}(t^*)=0$, after that the transfer reverses its direction
from 
$K_{\theta}$  to $K_{\nabla\phi}$
as $t\to \infty$, until it purges and gives $K_\theta(\infty)=0$
and $K(\infty)=K_{\nabla\phi}(\infty)$.

\end{document}